\documentclass[mathpazo]{cicp}

\usepackage{graphicx}
\usepackage{lipsum}
\usepackage{amsfonts}
\usepackage{graphicx}
\usepackage{epstopdf}
\usepackage{algorithmic}
\ifpdf
\DeclareGraphicsExtensions{.eps,.pdf,.png,.jpg}
\else
\DeclareGraphicsExtensions{.eps}
\fi

\usepackage{multicol}        
\usepackage[bottom]{footmisc}

\usepackage{amssymb}
\usepackage{amsmath}
\usepackage{amsthm}
\usepackage{wasysym}

\usepackage{bm}
\usepackage{xspace}

\usepackage{booktabs}
\usepackage{tabularx}

\usepackage{caption}
\usepackage{subfig}

\usepackage{mathptmx}
\usepackage{color}

\usepackage{hyperref}
\usepackage[capitalize]{cleveref}

\def\@thmcountersep{.}

\newtheorem{mydef}[theorem]{Definition}

\newtheorem{myremark}[theorem]{Remark}

\newcolumntype{L}[1]{>{\raggedright\arraybackslash}p{#1}} 
\newcolumntype{C}[1]{>{\centering\arraybackslash}p{#1}} 
\newcolumntype{R}[1]{>{\raggedleft\arraybackslash}p{#1}} 
\newcolumntype{Y}{>{\centering\arraybackslash}X} 
\newcolumntype{Z}{>{\raggedleft\arraybackslash}X} 

\newcommand{\centercell}[1]{\multicolumn{1}{c}{#1}}
\newcommand{\head}[1]{\centercell{#1}}

\newcommand{\wb}{well-balanced\xspace}
\newcommand{\nwb}{non-well-balanced\xspace}

\newcommand{\ud}{\textrm{d}}
\newcommand{\dd}[2]{\frac{\ud #1}{\ud #2}}
\newcommand{\df}[2]{\frac{\partial #1}{\partial #2}}
\newcommand{\con}{\bm{q}}

\newcommand{\imh}{{i - \frac{1}{2}}}
\newcommand{\iph}{{i + \frac{1}{2}}}
\newcommand{\jmh}{{j - \frac{1}{2}}}
\newcommand{\jph}{{j + \frac{1}{2}}}

\newcommand{\fl}{\bm{f}}
\newcommand{\gl}{\bm{g}}
\newcommand{\bw}{\bm{w}}

\newcommand{\bs}{\bm{s}}

\newcommand{\half}{\frac{1}{2}}
\newcommand{\const}{\textrm{const}}

\newcommand{\hfl}{\bm{F}}
\newcommand{\hgl}{\bm{G}}
\newcommand{\hbs}{\bm{S}}

\newcommand{\nhfl}{\hat{\hfl}}
\newcommand{\nhgl}{\hat{\hgl}}

\newcommand{\pe}{\bar{p}}

\newcommand{\rhoe}{\bar{\rho}}
\newcommand{\wrt}{w.r.t.\xspace}

\renewcommand{\vec}[1]{\mathbf{#1}}

\renewcommand{\epsilon}{\varepsilon}

\newcommand{\oo}[1]{\ensuremath{\frac{1}{#1}}}

\newcommand{\D}[2][ ]{\frac{\partial #1}{\partial #2}}

\newcommand{\DD}[2]{\frac{\partial #1}{\partial #2}}

\newcommand{\x}{\ensuremath{\vec x}\xspace}

\newcommand{\brunt}{Brunt--V\"ais\"al\"a\xspace}
\newcommand{\tBV}{\xspace\ensuremath{t_\text{BV}}\xspace}
\newcommand{\nablaad}{\xspace\ensuremath{\nabla_{\text{ad}}}\xspace}
\newcommand{\nablaext}{\xspace\ensuremath{\nabla_{\text{ext}}}\xspace}

\newcommand{\args}{\left(  t,\vec x\right)}

\newcommand{\hvF}{\hat{\vec F}}

\newcommand{\Lorho}{\ensuremath{\left\Vert\rho-\rho_0\right\Vert_{L^1}}}
\newcommand{\Lop}{\ensuremath{\left\Vert p-p_0\right\Vert_{L^1}}}
\newcommand{\Lov}{\ensuremath{\left\Vert v\right\Vert_{L^1}}}
\newcommand{\Lou}{\ensuremath{\left\Vert u\right\Vert_{L^1}}}

\newcommand\ie{i.e.\ }
\newcommand\eg{e.g.\ }

\numberwithin{theorem}{section}

\def\@thmcountersep{.}

\newcommand \subheadseparator {\midrule}
\newcommand{\subheadstyle}[1]{#1}
\makeatletter
\newcommand{\manuallabel}[2]{\def\@currentlabel{#2}\label{#1}}
\makeatother

\newcounter{subcount}
\newcommand{\usesubcount}{%
	\alph{subcount}%
	\addtocounter{subcount}{1}%
	
}

\newcommand\makesublabel[2]{\manuallabel{#1}{#2}(#2)}
\newcommand{\sublabel}[1]{\makesublabel{#1}{\usesubcount}}
\newcommand{\subtabref}[1]{\cref{tab:short_time}\ref{#1}}

\newcommand{\TheTitle}{Second order finite volume scheme for Euler equations with gravity which is well-balanced for general equations of state and grid systems}

\newcommand{\TheAuthors}{J.\ P.\ Berberich, P.\ Chandrashekar, C.\ Klingenberg, F.\ K.\  R\"opke}
\newcommand{\TheAuthorsFull}{Jonas P.~Berberich\affil{1}, Praveen Chandrashekar\affil{2}, Christian Klingenberg\corrauth\comma\affil{1}, and Friedrich K.\  R\"opke\affil{3}\comma\affil{4}}

\begin{document}
	\title{\TheTitle}
	
	\author[\TheAuthors]{\TheAuthorsFull}
	\address{	\affilnum{1} Dept.~of Mathematics, Univ.~of W\"urzburg, 97074 W\"urzburg, Germany\\
		\affilnum{2} TIFR Center for Applicable Mathematics, Bangalore 560065, India\\
		\affilnum{3} Zentrum f{\"u}r Astronomie der Universit{\"a}t Heidelberg, Institut f{\"u}r Theoretische Astrophysik, Philosophenweg 12, 69120 Heidelberg, Germany\\
		\affilnum{4}
		Heidelberger Institut f{\"u}r Theoretische Studien, Schloss-Wolfsbrunnenweg 35, 69118 Heidelberg, Germany
	}
	\email{{\tt klingen@mathematik.uni-wuerzburg.de} (C.\ Klingenberg)}

	\begin{abstract}
		We develop a second order well-balanced finite volume scheme for compressible Euler equations with a gravitational source term.
		The well-balanced property holds for arbitrary hydrostatic solutions of the corresponding Euler equations without any restriction on the equation of state. 
		The hydrostatic solution must be known a priori either as an analytical formula or as a discrete solution at the grid points.
		The scheme can be applied on curvilinear meshes and in combination with any consistent numerical flux function and time stepping routines. These properties are demonstrated on a range of numerical tests.
	\end{abstract}
	
	\ams{76M12, 
		65M08, 
		65M20, 
		35L65, 
		76N15, 
		76E20, 
		85-08. 
	}
	
	\maketitle

	\section{Introduction}
	The Euler equations with gravitational source term are used to model the flow of gases in different fields of physical sciences. Examples include weather prediction, climate modeling, and several astrophysical application such as the modeling of stellar interiors. In many of these applications, the gas state is close to the hydrostatic solution. Hydrostatic solutions are non-trivial steady state solutions of the Euler equations with gravity which can be described using the differential equation $\nabla p(\x)=-\rho(\x)\nabla(\phi(\x))$, where $p$ is the gas pressure, $\rho$ is the gas density and $\phi$ is a given external gravitational potential. In order to resolve the flow dynamics, which can be seen as small perturbations of the hydrostatic solution, one has to be able to maintain the corresponding hydrostatic solution with a sufficiently small error. Conventional methods introduce a significant discretization error when trying to compute small perturbations of the hydrostatic solution, especially on coarse meshes. Since the computational effort using sufficiently fine meshes can be too high, especially in three-dimensional simulations, special numerical techniques for this problem have been developed called well-balanced schemes. Well-balanced schemes are able to maintain hydrostatic solutions close to machine precision even on coarse meshes.
	
	Well-balanced schemes have been developed for the well-known shallow water equations with non-flat bottom topography. The equation describing steady state solutions in the shallow water equations is given in an explicit algebraic form which favors the development of well-balanced schemes. Some examples are \cite{Audusse2004,Noelle2007}.  There are also well-balanced schemes for related models, like \eg the Ripa model \cite{Desveaux2016b,Touma2015}. More recently, well-balanced schemes for Euler equations with gravitational source term have been developed. 
	This is more delicate than for shallow water equations since the hydrostatic solutions are given implicitly via a differential equation. For different equations of state (EoS) different hydrostatic solutions can be found. This led to the development of well-balanced schemes which are restricted to certain EoS and classes of hydrostatic solutions.
	Early work on this topic has been conducted by Cargo and Le Roux \cite{Cargo1994}. LeVeque and Bale \cite{Leveque1999} applied a quasi-steady wave-propagation algorithm on the Euler equations with gravitational source term to maintain isothermal hydrostatic solutions numerically. This method has been expanded to isentropic solutions in \cite{Leveque2011b}.
	Even for high order schemes well-balancing is necessary if solutions close to a hydrostatic solution are computed \cite{Xing2013}. A high order well-balanced scheme for isothermal hydrostatic solutions is introduced in \cite{Xing2013}. The scheme includes a modified weighted essentially non-oscillatory (WENO) reconstruction and a suitable way to discretize the source term. Based on this idea, a non-staggered central scheme for the same class of hydrostatic solutions has been proposed in \cite{Touma2016}. Compact reconstruction WENO methods are applied to achieve well-balancing in \cite{Ghosh2016}. Discontinuous Galerkin (DG) well-balanced methods have been developed in \cite{Li2016}, \cite{Li2018}, and \cite{Chandrashekar2017}. The well-balanced method proposed in \cite{Chertock2018} is based on a reformulation of the Euler equations with gravity discretized using a central scheme.
	
	Another approach to achieve well-balancing for Euler with gravitational source term is the development of well-balanced relaxation schemes, see \cite{Desveaux2014,Desveaux2016,Thomann2018} and references therein. Here stable approximate Riemann solvers are constructed for well-balancing. 
	Methods based on hydrostatic reconstruction were first developed by Audusse et.\ al.\ \cite{Audusse2004} for the shallow water equations. Later they have been adapted for the Euler equations with gravitational source term, see \eg \cite{Chandrashekar2015}. Early applications for weather prediction can be seen in \cite{Botta2004}.
	
	Most of the above mentioned well-balanced schemes for the Euler equations with gravitational source term are developed to balance a certain class of hydrostatic solutions, typically isothermal or polytropic (including isentropic) solutions. Many of these schemes also require the use of specific numerical flux functions. They can also only be successfully applied if a certain equation of state (EoS) is used, in most cases an isentropic or ideal gas EoS. Especially for astrophysical applications this limitation is too strict. Astrophysical EoS are very complex since they take into account a wide range of physical effects. Typically, there are no closed form expressions for these EoS. Hence, one cannot hope to calculate the hydrostatic solutions analytically.
	An approximately well-balanced scheme for general EoS has been developed in \cite{Kaeppeli2014} for the isentropic hydrostatic solution which is independent of the EoS. In \cite{Kaeppeli2016} K\"appeli and Mishra 
	generalized their scheme. It admits discrete hydrostatic solutions for arbitrary EoS, which are a second order approximation of the continuous hydrostatic solution.
	
	In physical and especially astrophysical applications certain properties are required. With the well-balanced method we present in this paper we can satisfy the following requirements:
	\begin{itemize}
		\item The method can be applied in the finite volume framework, which is suitable for solving the Euler equations due to their conservation property.
		\item The method can be used to balance any arbitrary hydrostatic solution satisfying any EoS. This includes EoS which are not given in closed form but in a table as is common for example in astrophysical simulations.
		\item The method can be combined with any numerical flux function.
		\item The method can be combined with any time-stepping routine including implicit ones.\footnote{By the term time-stepping routines we refer to ODE solvers which are applied to evolve the semi-discrete scheme obtained after spatial discretization by a finite volume method.}
		\item The method can be combined with any reconstruction routine.
		\item The method can be applied on arbitrary curvilinear grid geometries.
		\item The method can be implemented in a straight forward way which does not need restructuring of the code for classical finite volume codes.
	\end{itemize}
	This makes our method attractive for application, where for example specially designed numerical flux functions or other routines have to be used. Following ideas from \cite{Ghosh2015} and \cite{Chandrashekar2015}, the scheme is based on a hydrostatic reconstruction and a certain second order discretization of the gravitational source term.
	
	To apply our method, the considered hydrostatic solution has to be known a priori either as an analytical formula or as a discrete solution at the grid points. 
	In many relevant applications this is the case, since small dynamics on a hydrostatic background is examined in simulations. In these cases the underlying hydrostatic solution is known or can be integrated from the data numerically. 
	The scheme can be implemented in a manner that the source code can be used to study any hydrostatic solution. The user has to only supply a subroutine that gives the hydrostatic solution either as a function of the spatial coordinates or as discrete solution at the grid points. The latter can be realized by a routine reading discrete data from a table.
	Especially for astrophysical applications, where hydrostatic solutions to complex EoS appear, our scheme is useful. To resolve the small dynamics on the hydrostatic background it can be necessary to use so-called low Mach number numerical flux functions (e.g.\ \cite{Miczek2015,Barsukow2017}). These can be used in combination with our well-balanced method to obtain a finite volume scheme which is both well-balanced and capable of resolving low Mach number fluxes. An example of such an application can be a simulation of a star with convective and non-convective regions. Convective stability of hydrostatic solutions is discussed in \cref{sec:brunt}. The capability of our method to correctly represent convective stability is demonstrated in \cref{sec:numerical}.
	
	The rest of the paper is structured as follows: In \cref{sec:EulerEquations} the Euler equations with gravitational source term are introduced for curvilinear coordinates. A few examples of EoS are presented in~\cref{sec:eos}. Hydrostatic solutions of the Euler equations with gravitational source term are discussed in \cref{sec:hystat} in general. Particular hydrostatic solutions are presented together with a method to determine stability properties of several equilibria. In \cref{sec:wbScheme} the well-balanced scheme is introduced and the well-balanced property is shown analytically. Various numerical experiments confirming the well-balanced property follow in \cref{sec:numerical}. Also, the ability of the scheme to resolve small perturbations on the hydrostatic solution is shown using numerical tests and the accuracy of the scheme is tested. Finally, conclusions of the paper are drawn in \cref{sec:conclusions}.

	\section{2D Euler Equations with Gravitational Source Term}
	\label{sec:EulerEquations}
	The 2D Euler equations which model the balance laws of mass, momentum, and energy under the influence of gravity in Cartesian coordinates are given by
	\begin{equation}
	\df{\con}{t} + \df{\fl}{x} + \df{\gl}{y} = \bs,
	\label{eq:eul2d_cartesian}
	\end{equation}
	where the conserved variables, fluxes and source terms are
	\[
	\con = \begin{bmatrix}
	\rho \\ \rho u \\ \rho v \\ E \end{bmatrix}, \quad
	\fl = \begin{bmatrix}
	\rho u \\
	p + \rho u^2 \\
	\rho u v \\
	(E+p)u \end{bmatrix}, \quad 
	\gl = \begin{bmatrix}
	\rho v \\
	\rho u v \\
	p + \rho v^2 \\
	(E+p)v \end{bmatrix}, \quad
	\bs = \begin{bmatrix}
	0 \\
	-\rho \df{\phi}{x} \\
	-\rho \df{\phi}{y} \\
	0
	\end{bmatrix}
	\]
	wit $\rho,p>0$.
	Moreover, $E=\rho\varepsilon + \tfrac{1}{2}\rho|\vec v|^2 + \rho\phi$ is the total energy per unit volume with $\vec v=(u,v)^T$ being the velocity and $\epsilon$ the specific internal energy. The scalar function $\phi$ is a given gravitational potential. If the computational domain has curved boundaries, we can use body-fitted coordinates $\boldsymbol\xi=(\xi,\eta)$ which we refer to as curvilinear coordinates. We map the physical domain in $\x$-space by a smooth mapping to a unit square in the $\boldsymbol\xi$-space in a bijective way. The Euler equations can be transformed to curvilinear coordinates and written as
	\begin{equation}
	\df{}{t}(J\con) + \df{}{\xi}[(x_\eta^2+y_\eta^2)^{\tfrac{1}{2}} \hfl] + \df{}{\eta}[(x_\xi^2+y_\xi^2)^{\tfrac{1}{2}} \hgl] = \hbs,
	\label{eq:eul2d}
	\end{equation}
	where
	$\hfl = \ell_1 \fl + \ell_2 \gl$,
	$\hgl = m_1 \fl + m_2 \gl$,
	$\hbs = J \bs$.
	
	For details see \cref{appendix:trafo}.
	The subscripts $(\cdot)_\xi$, $(\cdot)_\eta$ denote partial derivatives. 
	In the above equations $(\ell_1, \ell_2)$ and $(m_1,m_2)$ are unit vectors defined as
	$(\ell_1, \ell_2) = (y_\eta, -x_\eta)(x_\eta^2+y_\eta^2)^{-\half}$ and $(m_1,m_2) = (-y_\xi, x_\xi)(x_\xi^2+y_\xi^2)^{-\half}$.
	$J$ is the Jacobian of the transformation $(\xi,\eta) \to (x,y)$ given by
	$J = x_\xi y_\eta - x_\eta y_\xi$.
	The transformed fluxes can be written as
	\[
	\hfl =  \begin{bmatrix}
	\rho U \\
	p \ell_1 + \rho u U \\
	p \ell_2 + \rho v U \\
	(E+p)U
	\end{bmatrix}, \qquad \hgl = \begin{bmatrix}
	\rho V \\
	p m_1 + \rho u V \\
	p m_2 + \rho v V \\
	(E+p)V
	\end{bmatrix},
	\]
	where
	$U = u \ell_1  + v \ell_2$ and $V = u m_1 + v m_2$
	
	\subsection{Equations of State}
	\label{sec:eos}
	
	In the Euler equations \cref{eq:eul2d}, there are five independent scalar unknowns: $\rho$, $\rho u$, $\rho v$, $\rho E$, $p$. Yet, the system consists of only four scalar equations. In order to close the system, we need one more equation relating at least two of the above quantities. The common choice to close the system is using an EoS. 
	\begin{mydef}
		A relation between the quantities $\rho$, $p$, and $\epsilon$, given in an explicit or implicit form, is called \emph{equation of state} (EoS).
		\label{def:eos}
	\end{mydef}
	
	We proceed with giving two examples of EoS. 
	
	\subsubsection{Ideal Gas}
	\label{sec:ideal_gas}
	
	A widely applicable EoS derives from considering an ideal classical gas. The pressure for the ideal gas EoS is given by
	\begin{equation}
	p(\rho,\epsilon) := p(\rho,T(\rho,\epsilon)) := \frac{R}{\mu}\rho T(\rho,\epsilon),
	\label{eq:ideal_eos}
	\quad\text{where}\quad
	T(\rho,\epsilon):=\frac{(\gamma-1)\mu}{R}\epsilon.
	\end{equation}
	We set the gas constant $R$ and the mean molecular weight $\mu$ to $R=\mu=1$. For the tests presented in \cref{sec:numerical}, the specific heat ratio $\gamma$ is set to $\gamma=1.4$, which is suitable to describe a diatomic gas (\eg air). For this paper, it is convenient to write the ideal gas EoS in the form \cref{eq:ideal_eos} depending on the \emph{temperature} $T$ instead of the explicit dependence of $\epsilon$. We use this form of the EoS to formulate hydrostatic solutions with certain temperature profiles in \cref{sec:numerical}.
	
	\subsubsection{Ideal Gas with Radiation Pressure}
	\label{sec:radiative_eos}
	
	At very high temperatures, such as in stellar interiors, radiation pressure can be relevant. For this case, we add a term corresponding to the Stefan--Boltzmann law to the gas pressure given in \cref{eq:ideal_eos}. Again, let us set all physical constants except $\gamma$ to such values that all prefactors vanish. The EoS for an ideal gas with radiation pressure is then given by
	\begin{equation}
	p(\rho,\epsilon) := p(\rho,T(\rho,\epsilon)):= \rho T(\rho,\epsilon) + T^4(\rho,\epsilon)
	\label{eq:radiative_eos}
	\end{equation}
	(\eg \cite{Chandrasekhar1958}),
	where the temperature $T$ is defined implicitly via
	\begin{equation}
	\epsilon = \frac{T(\rho,\epsilon)}{\gamma-1} + \frac{3}{\rho} T^4(\rho,\epsilon).
	\label{eq:radiative_epsilon}
	\end{equation}
	The relation \cref{eq:radiative_epsilon} cannot be rewritten to express $T$ explicitly. It has to be solved numerically, \eg using Newton's method. The radiation constant $a$ that is present in the source \cite{Chandrasekhar1958} is set to $3$ to obtain these equations.
	
	\section{Hydrostatic Solutions}
	\label{sec:hystat}
	A hydrostatic solution of \cref{eq:eul2d} is a solution in which the fluid is at rest and the density and pressure are independent of time; we will denote the hydrostatic solution using $\bar{(\cdot)}$. Since the velocity is zero, the continuity and energy equations are automatically satisfied, and the momentum equation reduces to
	\begin{equation}
	\nabla \bar p = -\bar\rho\nabla\phi.
	\label{eq:hystat}
	\end{equation}
	This relation is called \emph{hydrostatic equation}. In order to solve this equation we will assume an EoS to relate pressure and density. Moreover, since an equation of state can depend on the internal energy, we need to make additional assumptions in order to be able to solve the above hydrostatic equation. Examples of such assumptions are given in \cref{sec:hystat_examples}.
	
	We intent to test our well-balanced numerical method on stable hydostatic equilibria. Thus in the next subsection we introduce the notion of stability we shall be using.
	
	\subsection{\brunt Frequency}
	\label{sec:brunt}
	
	Consider a hydrostatic solution of the Euler equations with a gravitational source term given via functions $\rho:=\rho_\text{ext}(\chi)$ and $p:=p_\text{ext}(\chi)$. Hypothetically, we take a fluid element from a position $\chi=\chi_0$ and displace it adiabatically, \ie without an exchange of heat with the surrounding fluid. According to \cite{Maeder2009} the equation of motion for a fluid element is
	\begin{equation*}
	\rho_\text{int} \frac{d^2\chi}{dt^2} + g \left.\left( \D[\rho_\text{int}]{\chi}-\rho_\text{ext}^\prime(\chi) \right)\right|_{\chi_0}(\chi-\chi_0)=0
	\end{equation*}
	for a small displacement $(\chi-\chi_0)$. $\rho_\text{int}$ is the density of the displaced fluid element, $g=\phi^\prime(\chi)$ is the gravity in $\chi$-direction.
	The solution of this ordinary differential equation is of the form
	\begin{equation}
	\chi-\chi_0 = a\,e^{iNt}
	\label{eq:brunt_sol}
	\end{equation}
	with the initial $\chi$-displacement being equal to $a$ and the \emph{\brunt frequency} $N$ given by
	\begin{equation}
	N^2 = \left.\frac g{\rho_\text{ext}}\left( \D[\rho_\text{int}]{\chi}-\rho_\text{ext}^\prime(\chi) \right)\right|_{\chi_0}.
	\label{eq:brunt}
	\end{equation}
	A form which is more practical in some situations is derived e.g.~ in \cite{Maeder2009}:
	\begin{equation}
	N^2 = \frac g{H_p}\left( \nablaad - \nablaext \right).
	\label{eq:brunt_ad}
	\end{equation}
	The \emph{pressure height scale} is defined as
	\begin{equation}
	H_p := -\frac{p_\text{ext}(\chi)}{p_\text{ext}^\prime(\chi)}=\frac{p_\text{ext}(\chi)}{\rho_\text{ext}\phi^\prime(\chi)}
	\label{def:Hp}
	\end{equation}
	for a hydrostatic solution and the temperature gradients are
	\begin{equation}
	\nablaad:=\left.\frac{p_\text{int}}{T_\text{int}}\frac{\partial T_\text{int}}{\partial  p_{\text{int}}}\right|_{s}, \qquad \nablaext:=\frac{p_\text{ext}T_\text{ext}^\prime(\chi)}{T_\text{ext}p_\text{ext}^\prime(\chi)}.
	\label{def:nablas}
	\end{equation}
	The derivative in the adiabatic temperature gradiant $\nablaad$ is taken at constant entropy.
	Note that $\nablaad$ depends on the EoS, while $\nablaext$ depends on the stratification of the data.
	
	Assuming an ideal gas EoS, one can find (e.g.~\cite{Maeder2009})  
	\begin{equation}
	\nablaad=\frac{\gamma-1}{\gamma}.
	\label{eq:ad_ideal}
	\end{equation}
	According to \cite{Kippenhahn1990} (their Eqs. (13.2), (13.3), and (13.12)), for an ideal gas with radiation pressure the adiabatic temperature gradient can be given as
	\begin{equation}
	\nablaad=\frac{2 p (4 p-3 \rho  T)}{32 p^2-24 p \rho  T-3 \rho ^2 T^2}.
	\label{eq:ad_radiative}
	\end{equation}

	\begin{mydef}
		\label{def:convective_stability}
		We call a hydrostatic solution of the Euler equations with gravitational source term \emph{stable with respect to convection}, if $N^2\geq0$ on the whole domain $\Omega$. Otherwise we call it \emph{unstable with respect to convection}.
	\end{mydef}
	If $N^2>0$, $N$ is real and \cref{eq:brunt_sol} describes an oscillation around $\chi_0$ with the frequency $N$. If $N=0$, the fluid element does not move at all. In both cases, the amplitude never surpasses the initial displacement $a$. If $N^2<0$, $N$ is imaginary, so that the exponent in \cref{eq:brunt_sol} is $\pm i N t$. $(\chi-\chi_0)$ is then a real exponential function with a positive exponent in one of the solutions. The displacement increases in time in that case.
	In some numerical tests in in which we consider hydrostatic solutions which are stable \wrt convection we use the minimal \brunt time
	\begin{equation}
	t_\text{BV}:=\min_{\x \in\Omega} t_\text{BV}^\text{loc}(\x):=\frac{2\pi}{\max_{\x\in\Omega} N(\x)}
	\label{eq:brunt_time}
	\end{equation}
	as a reference time. It seems to be a natural time scale for the evolution of small perturbations, such as numerical errors, on a hydrostatic solution. The \brunt frequency and time depend on the spatial position $\x$, since $\rho$ and $\phi$ depend on $\x$. We need a global time as reference time. The minimal \brunt time is used because it corresponds to the fastest oscillations. One can expect that significant spurious dynamics, introduced via small statistical perturbations, first occur around the coordinates with the highest value for $N$ and lowest value for $t_\text{BV}^\text{loc}$ if $N^2>0$.
	
	\begin{myremark}
		The definition of the \brunt frequency can be extended to hydrostatic solutions with genuinely multidimensional smooth $\rho$, $p$, and $\phi$ by locally aligning the $\chi$-coordinate with $\nabla \phi$.
		\label{remark:brunt}
	\end{myremark}
	
	\subsection{Examples of Hydrostatic Solutions considering Stability w.r.t. Convection}
	\label{sec:hystat_examples}
	\subsubsection{Polytropic and Isentropic Solutions}
	\label{sec:polytropic}
	Polytropic solutions of \cref{eq:hystat} are of the form
	\begin{equation}
	\theta(\x):=1-\frac{\nu-1}{\nu}\phi(\vec x),\quad\bar\rho(\vec x) = \theta(\x)^{\frac{1}{\nu-1}},\quad\bar p(\vec x)=\bar\rho(\x)^\nu
	\label{eq:polytropic}
	\end{equation}
	with $\nu>0$. \Cref{eq:polytropic} is a solution of \cref{eq:hystat} to several different EoS, since there is freedom for the choice of the internal energy stratification.
	Using an ideal gas EoS and $\nu=\gamma$ defines the special case of an isentropic solution. In the case of an ideal gas EoS, we have $\bar T\equiv\theta$ for the equilibrium temperature $\bar T$.
	
	The stability with respect to convection can be computed analytically in the case of an ideal gas EoS using \cref{eq:brunt_ad,eq:ad_ideal}. 
	For \cref{eq:polytropic} with an ideal gas EoS we have
	\begin{equation*}
	\frac{\bar p(\x)}{\bar \rho(\x)}=\bar T(\x)=1-\frac{\nu-1}{\nu}\phi(\x)
	\end{equation*}
	and
	\begin{equation*}
	\ln \bar p = \ln\left( \bar T^{\frac\nu{\nu-1}} \right) = \frac\nu{\nu-1}\ln \bar T\quad\Rightarrow\quad\nablaext=\frac{\nu-1}\nu.
	\end{equation*}
	Using \cref{eq:brunt_ad,eq:ad_ideal} we get
	\begin{equation}
	N(\x)^2 = \left( 1-\frac{\nu-1}{\nu}\phi(\x) \right)^{-1}\left( \frac{\gamma-1}{\gamma}-\frac{\nu-1}{\nu} \right). 
	\label{eq:N_polytropic}
	\end{equation}
	Assuming $\nu\geq 1$ and $\phi\leq1<\frac{\nu}{\nu-1}$, the polytropic solution is stable with respect to convection if and only if $\nu\leq\gamma$. In the isentropic case, the equilibrium is marginally stable everywhere, since $N^2\equiv0$. An adiabatically displaced fluid element neither oscillates around the original position nor moves away from it, but it just holds the new position. A \brunt time can not be calculated, since $N=0$ is the denominator in \cref{eq:brunt_time}.
	
	In the case $\nu<\gamma$, the criterion in \cref{def:convective_stability} is satisfied. The global \brunt time, calculated using \cref{eq:brunt_time} and \cref{eq:N_polytropic}, is
	\begin{equation*}
	\tBV=\max_{\x\in\Omega}\left( \sqrt{\frac{4\pi^2\left( 1-\frac{\nu-1}{\nu}\phi(\x) \right)}{\left( \frac{\gamma-1}{\gamma}-\frac{\nu-1}{\nu} \right)}} \right)
	= 
	\frac{2\pi}{\sqrt{\left( \frac{\gamma-1}{\gamma}-\frac{\nu-1}{\nu} \right)}}. 
	\end{equation*}
	In the case $\nu>\gamma$ no \brunt time can be given since the solution is then unstable \wrt convection.
	
	If the ideal gas with radiation pressure EoS is used, stability w.r.t.\ convection can be discussed using \cref{eq:brunt_ad} together with \cref{eq:ad_radiative} on the data. It can not be done analytically since there is no explicit expression for the temperature in this case. Instead, $N^2$ can be approximated numerically from the data on the grid. For the spatial derivative of $T$ we will use a central finite difference approximation.
	
	\subsubsection{Isothermal Solutions}
	\label{sec:isothermal}
	Hydrostatic solutions with constant temperature, in our case $\bar T\equiv 1$, are called isothermal solutions. An isothermal solution for the ideal gas EoS \cref{eq:ideal_eos} is given by
	\begin{equation}
	\bar \rho(\x) = \bar p(\x) = \exp(-\phi(\x)).
	\label{eq:isothermal}
	\end{equation} 
	One can also find an isothermal solution for the ideal gas with radiative pressure EoS \cref{eq:radiative_eos}. It is given by
	\begin{equation}
	\bar \rho(\x) = \exp\left( -\phi(\x) \right),\quad \bar p =  \exp\left( -\phi(\x) \right) + 1.
	\label{eq:radiative_isothermal}
	\end{equation}
	Again, we can analyze the stability with respect to convection in the case of an ideal gas EoS. Since $T\equiv \const$, we have $\rho\equiv p$ and $\nablaext\equiv0$ for this solution. With that, \cref{eq:brunt_ad,eq:ad_ideal} yield a spatially constant \brunt frequency of
	\begin{equation*}
	N=\sqrt{\frac{\gamma-1}{\gamma}},\quad
	\text{and hence}\quad
	\tBV = 2\pi\sqrt{\frac{\gamma}{\gamma-1}},
	\end{equation*}
	which is an implication of \cref{eq:brunt_time}. This hydrostatic solution is always stable with respect to convection.
	
	Plugging the hydrostatic solution into \cref{eq:brunt_ad,eq:ad_radiative} we find
	\begin{equation}
	N^2(\x) = 
	\frac{3}{4 \left(8 \exp(\phi(\x)) \left(4 \exp(\phi(\x) )+5\right)+5\right)}+\frac{1}{4},
	\label{eq:brunt_rad_isothermal}
	\end{equation}
	which is positive for all $\x$. The isothermal hydrostatic solution for the ideal gas with radiation pressure EoS is hence also stable w.r.t.\ convection.

	\subsubsection{The tanh-Profile}
	\label{sec:tanhProfile}
	
	To demonstrate the flexibility of our scheme in \cref{sec:2d_tanh}, we introduce a hydrostatic solution which can for example be found in \cite{Edelmann2014}.
	We assume a temperature profile
	\begin{equation*}
	\bar T(\x)=\bar T(x,y)=1+\Delta T \tanh\left( \frac{x}{\mu} \right)
	\end{equation*}
	with constants $\Delta T,\mu > 0$. For this profile the following solution of the hydrostatic equation \cref{eq:hystat} together with the ideal gas EoS \cref{eq:ideal_eos} can be found:
	\begin{equation}
	\bar p(\x)=\bar p(x,y)=\exp\left( -\frac{x-\Delta T\mu\log\left(\cosh\left(\frac{x}{\mu}\right)+\Delta T\sinh\left( \frac{x}{\mu} \right)\right)}{1-\Delta T^2} \right),\;  \bar\rho(\x)=\frac{\bar p(\x)}{\bar T(\x)}.
	\label{eq:tanh_rhop}
	\end{equation}
	The corresponding gravitational potential is $\phi(\x)=\phi(x,y)=x$.
	The ideal gas EoS \cref{eq:ideal_eos} is used, hence the \brunt frequency can be computed via \cref{eq:brunt_ad,eq:ad_ideal} if $\nablaext$ can be determined.
	We can approximate the maximal $N^2$ numerically using the values on the initial numerical grid and numerical derivatives.
	
	\subsubsection{Numerically Approximated Hydrostatic Solution}
	\label{sec:integratedProfile}
	
	Hydrostatic solutions which are analytically known are a special case. In astrophysical applications it is common, that underlying hydrostatic solutions are not known analytically. Instead, initial data for the simulation are given in the form of data points on the grid. Knowing, that these data are close to a hydrostatic solution a close by hydrostatic solution can be obtained by numerically solving the hydrostatic equation and EoS. For that one can for example use the given temperature profile and solve the system for density and pressure. 
	
\begin{figure}
		\centering
		\includegraphics[scale=1.]{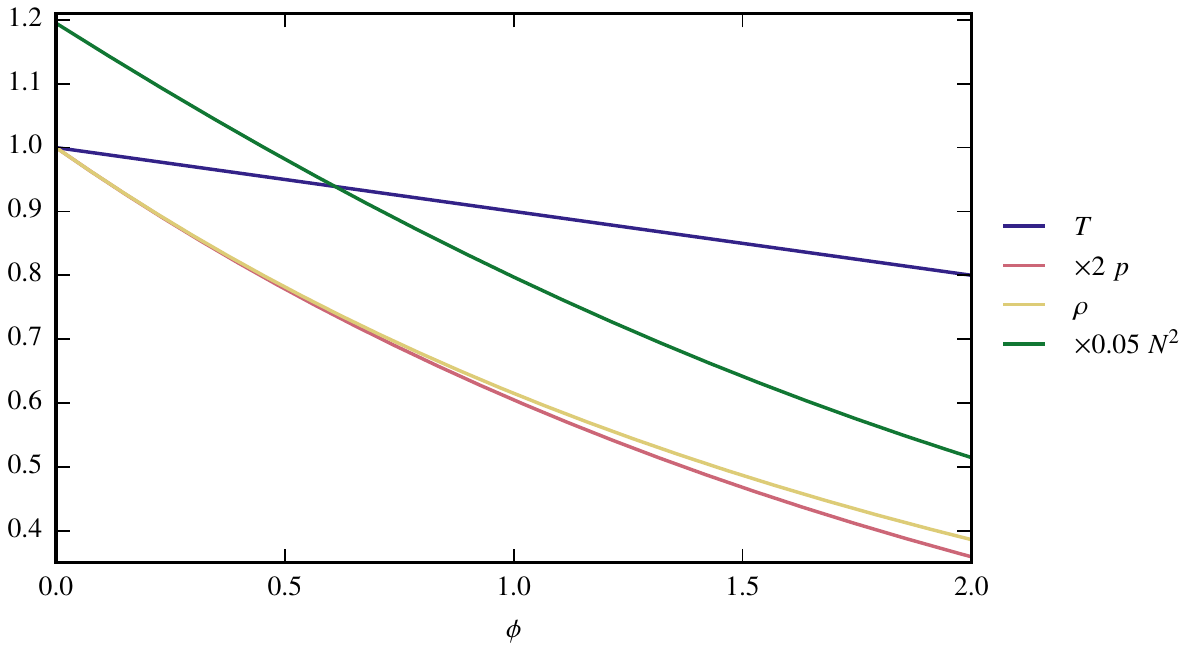}
		\caption{\label{fig:integrated_data} Data of the integrated hydrostatic solution presented in \cref{sec:integratedProfile}.}
	\end{figure}

	As an example we assume the following given data: Let the gravitational potential be given as $\phi(\x)=\phi(x,y)=x+y$ and the temperature $T(\x)=1-0.1\phi(\x)$. We assume an ideal gas with radiation pressure as given in \cref{eq:radiative_eos,eq:radiative_epsilon}. Using Chebfun \cite{Driscoll2014} in the numerical software MATLAB we solve the hydrostatic equation and EoS for density and pressure corresponding to the given temperature profile. The result is shown in \cref{fig:integrated_data}. 
	
	To discuss stability w.r.t.\ convection we use the criterion given by \cref{eq:brunt_ad,eq:ad_radiative}.
	In our case, the direction $\chi$ is the diagonal direction $\chi(x,y)=\frac1{\sqrt{2}}(x+y)$. The resulting values for $N^2$ are shown in \cref{fig:integrated_data}. We see that the hydrostatic solution we computed is stable w.r.t.\ convection, since the criterion $N^2>0$ is met in the whole domain.
	\section{The $\alpha$-$\beta$ Well-Balanced Scheme}
	\label{sec:wbScheme}
	In the following we construct our well-balanced scheme for Euler equations with gravitational source term on curvilinear grids.
	
	\subsection{Hydrostatic Solution and Source Terms}
	
	Let us rewrite the hydrostatic solution as
	\[
	\rhoe(x,y) = \rho_0 \alpha(x,y), \qquad \pe(x,y) = p_0 \beta(x,y),
	\]
	which is merely a rescaling of the hydrostatic solution so that $\alpha$, $\beta$ are positive dimension-less quantities, $\rho_0$ and $p_0$ are positive constants. The non-dimensional functions will be useful later during the reconstruction step of the finite volume method. Using \cref{eq:hystat}, the gravitational force can be written as
	\[
	\df{\phi}{x} = - \frac{p_0}{\rho_0 \beta} \df{\alpha}{x}, \qquad \df{\phi}{y} = - \frac{p_0}{\rho_0 \beta} \df{\alpha}{y}.
	\]
	We will now rewrite the source terms in the general coordinates $(\xi,\eta)$ as follows:
	\begin{eqnarray*}
		-J \rho \df{\phi}{x} &=&  \frac{p_0 \rho}{\beta} J \df{\alpha}{x} 
		=  \frac{p_0 \rho}{\beta} J (\alpha_\xi \xi_x + \alpha_\eta \eta_x)
		\stackrel{\cref{eq:trans}}= \frac{p_0 \rho}{\beta} (\alpha_\xi y_\eta - \alpha_\eta y_\xi ) \\
		&=& \frac{p_0 \rho}{\rho_0 \beta}( (\alpha y_\eta)_\xi - (\alpha y_\xi)_\eta),
	\end{eqnarray*}
	where the last relation holds due to $(y_\eta)_\xi = (y_\xi)_\eta$. In the same way we get
	\[
	-J \rho \df{\phi}{y} = \frac{p_0 \rho}{\rho_0 \beta} (- (\alpha x_\eta)_\xi + (\alpha x_\xi)_\eta ).
	\]
	This re-formulation of the source term is one of the key elements in achieving a well-balanced scheme.
	\subsection{Finite Volume Scheme}
	\label{sec:wbScheme_fv}
	The body-fitted curvilinear coordinates define a mesh for the physical domain. The grid points are indexed by $(i,j)$ and  the $(i,j)$'th cell is denoted $C_{i,j}$ in the physical space and by $\hat{C}_{i,j}$ in the computational space, i.e., in $(\xi,\eta)$ coordinates.  The semi-discrete finite volume scheme is obtained by integrating \cref{eq:eul2d} over the cell $\hat{C}_{i,j}$ leading to
	\begin{equation}
	\Omega_{i,j} \dd{}{t}\con_{i,j} + \nhfl_{\iph,j} A_{\iph,j} - \nhfl_{\imh,j} A_{\imh,j} + \nhgl_{i,\jph} A_{i,\jph} - \nhgl_{i,\jmh} A_{i,\jmh} = \hbs_{i,j},
	\label{eq:fvm2d}
	\end{equation}
	where
	\[
	\Omega_{i,j} = \int_{\hat{C}_{i,j}} J \ud\xi \ud\eta = \int_{C_{i,j}} \ud x \ud y 
	\]
	is the physical area of the cell. The numerical fluxes $\nhfl_{\iph,j} = \nhfl(\con_{\iph,j}^L,\con_{\iph,j}^R, \ell_{\iph,j})$, $\nhgl_{i,\jph} = \nhgl(\con_{i,\jph}^L,\con_{i,\jph}^R, m_{i,\jph})$, etc. are consistent with the fluxes $\hfl$, $\hgl$, and the unit normal vectors are given by
	\[
	\ell_{\iph,j} = \frac{1}{A_{\iph,j}} \left[ (y_{\iph,\jph} - y_{\iph,\jmh}), \ -(x_{\iph,\jph} - x_{\iph,\jmh}) \right],
	\]
	\[
	m_{i,\jph} = \frac{1}{A_{i,\jph}} \left[ -(y_{\iph,\jph} - y_{\imh,\jph}), \ (x_{\iph,\jph} - x_{\imh,\jph}) \right],
	\]
	with
	\[
	A_{\iph,j} = \sqrt{ (x_{\iph,\jmh} - x_{\iph,\jph})^2 + (y_{\iph,\jmh} - y_{\iph,\jph})^2},
	\]
	\[
	A_{i,\jph} = \sqrt{ (x_{\iph,\jph} - x_{\imh,\jph})^2 + (y_{\iph,\jph} - y_{\imh,\jph})^2}.
	\]
	The gravitational source term is discretized as
	\begin{eqnarray*}
		S_{i,j}^{(1)} &=& 0, \\
		S_{i,j}^{(2)} &=&  +\frac{p_0 \rho_{i,j}}{\rho_0 \alpha_{i,j}} \left[ (y_\eta)_{\iph,j} \beta_{\iph,j} - (y_\eta)_{\imh,j} \beta_{\imh,j} \right] \\&& -  \frac{p_0 \rho_{i,j}}{\rho_0 \alpha_{i,j}} \left[ (y_\xi)_{i,\jph} \beta_{i,\jph} - (y_\xi)_{i,\jmh} \beta_{i,\jmh} \right],\\
		S_{i,j}^{(3)} &=& - \frac{p_0 \rho_{i,j}}{\rho_0 \alpha_{i,j}}  \left[ (x_\eta)_{\iph,j} \beta_{\iph,j} - (x_\eta)_{\imh,j} \beta_{\imh,j} \right] \\&& + \frac{p_0 \rho_{i,j}}{\rho_0 \alpha_{i,j}}  \left[ (x_\xi)_{i,\jph} \beta_{i,\jph} - (x_\xi)_{i,\jmh} \beta_{i,\jmh} \right], \\
		S_{i,j}^{(4)} &=& 0.
	\end{eqnarray*}
	The derivatives of the mesh in the source term are computed using central differences, e.g.
	\[
	(y_\eta)_{\iph,j} = y_{\iph,\jph} - y_{\iph,\jmh}, \quad \textrm{etc.},
	\]
	where
	\[
	y_{\iph,\jph} = \frac{1}{4}( y_{i,j} + y_{i+1,j} + y_{i+1,j+1} + y_{i,j+1}), \quad \textrm{etc.}
	\]
	To obtain the values at the face $\con_{\iph,j}^L$, $\con_{\iph,j}^R$, $\con_{i,\jph}^L$, $\con_{i,\jph}^R$ used to compute the numerical flux, we reconstruct the following set of variables
	\begin{equation}
	\label{eq:hystat_variables}
	\bw = [ \rho/\alpha, \ u, \ v, \ p/\beta ]^\top.
	\end{equation}
	We obtain the reconstructed primitive variables $\vec u$ via
	\begin{equation}
	\label{eq:transform_reconstructed}
	\vec u_{i+\half,j}^{L,R} = T\left(\vec x_{i+\half,j}\right)\bw_{i+\half,j}^{L,R}, \qquad
	\vec u_{i,j+\half}^{L,R} = T\left(\vec x_{i,j+\half}\right)\bw_{i,j+\half}^{L,R},
	\end{equation}
	where
	$
	T(\vec x) = \text{diag} \left(\alpha(\vec x)^{-1}, 1, 1, \beta(\vec x)^{-1} \right)
	$
	and transform them to conservative variables in the canonical way if necessary in the implementation.
	\subsection{Well-Balanced Property}
	The following theorem establishes the well-balanced property of the 2D finite volume scheme.
	
	\begin{theorem}
		The finite volume scheme (\cref{eq:fvm2d}) together with any consistent numerical flux and reconstruction of $\bw$ variables is well-balanced in the sense that the initial condition given by
		\begin{equation}
		u_{i,j} = v_{i,j} = 0, \qquad \rho_{i,j}/\alpha_{i,j} = \rho_0 = \const, \qquad p_{i,j}/\beta_{i,j} = p_0 = \const, \qquad \forall \ (i,j)
		\label{eq:wb2d}
		\end{equation}
		is preserved by the numerical scheme.
	\end{theorem}
	
	\noindent
	\underline{Proof}: Let us start the computations with an initial condition which is hydrostatic and hence satisfies~(\cref{eq:wb2d}). Since the $\bw$ variables are constant for this initial condition, any reconstruction scheme will be exact and give
	\[
	\rho_{\iph,j}^L = \rho_{\iph,j}^R = \rho_{\iph,j}, \quad p_{\iph,j}^L = p_{\iph,j}^R = p_{\iph,j},
	\]
	\[
	u_{\iph,j}^L = u_{\iph,j}^R = v_{\iph,j}^L = v_{\iph,j}^R = 0, 
	\]
	etc., where $\rho_{\iph,j} = \rho_0 \alpha_{\iph,j}$, $p_{\iph,j} = p_0 \beta_{\iph,j}$, etc. By consistency of the flux we get
	\begin{align*}
	\nhfl_{\iph,j} A_{\iph,j} &= \begin{bmatrix}
	0 \\
	+(y_{\iph,\jph}-y_{\iph,\jmh}) p_{\iph,j} \\
	-(x_{\iph,\jph}-x_{\iph,\jmh}) p_{\iph,j} \\
	0 \end{bmatrix}, \\ 
	\hat{\vec G}_{i,\jph} A_{i,\jph} &= \begin{bmatrix}
	0 \\
	-(y_{\iph,\jph}-y_{\imh,\jph}) p_{i,\jph} \\
	+(x_{\iph,\jph}-x_{\imh,\jph}) p_{i,\jph} \\
	0 \end{bmatrix}.
	\end{align*}
	An inspection of the source terms shows that they are exactly equal to the flux terms, leading to a well-balanced scheme.
	
	\begin{myremark}
		The well-balanced scheme presented in this section allows the choice of arbitrary consistent methods for the numerical flux function, the reconstruction method, the time-stepping routine, and the grid topology.  Note that the source term discretization is a second order approximation. The particular form of the source term discretization is of central importance for our well-balanced method. Also, in Eq. (4.1) the fluxes are evaluated as interface-averages. This limits the resulting scheme to second order accuracy.
	\end{myremark}
	\begin{myremark}
		The proposed well-balanced scheme excels due to generality. The freedom to write any hydrostatic solution into the functions $\alpha$ and $\beta$ makes the scheme flexible and applicable for all EoS. Only the choice of the reconstructed variables and the source term discretization have to be adapted to implement this scheme. This leaves freedom in the choice of all other components such as time-stepping, numerical flux function, reconstruction method, and grid geometry.
	\end{myremark}

	\subsection{Numerical Flux Function}
	
	We use Roe's approximate Riemann solver \cite{Roe1981} for the numerical flux function which is given by
	\begin{align*}
	\hvF_\chi(\vec q_L,\vec q_R) &= \half\left[ \vec F_\chi(\vec q_L) + \vec F_\chi(\vec q_R) - D_{\chi,\text{Roe}} (\vec q_R-\vec q_L)\right]
	\label{def:RoeTypeFlux}
	\end{align*}
	for $\chi=x,y$ and the fluxes $\vec F_x=\vec F, \vec F_y = \vec G$.
	The matrix $D_{\chi,\text{Roe}}$ is the Roe diffusion matrix. 
	With the eigenvalues $\lambda_\chi^l$ and eigenvectors $\vec r_\chi^l$ for $l=1,\dots,4$ of the flux Jacobian $A_\chi=\DD{\vec F_\chi}{\vec q}$ we can define the matrix
	\begin{equation}
	D_{\chi,\text{Roe}}:=|A_\chi|:= R_\chi|\Lambda_\chi|R_\chi^{-1},\qquad \text{where }|\Lambda_\chi|:=\text{diag}(|\lambda_\chi^l|),\quad R_\chi:=(\vec r_\chi^1,\dots,\vec r_\chi^4).
	\label{}
	\end{equation}
	The Diffusion matrix $D_{\chi,\text{Roe}}=D_{\chi,\text{Roe}}(\vec q^\prime)$ depends on a state vector $\vec q^\prime$. It is evaluated at the Roe-average state $\vec q^\prime=\vec q_\text{Roe}(\vec q^L,\vec q^R)$, which is given in \cite{Roe1981}.

	\subsection{Reconstruction}
	
	We use two types of reconstruction schemes in this paper, but this is not the only possible choice as the scheme is well-balanced for any reconstruction scheme. The constant reconstruction scheme, applied on the  variables $\bw$, is given by 
	\begin{align*}
	\bw^L_{i+\half,j} &= \bw_{i,j},& \bw^R_{i+\half,j} &= \bw_{i+1,j},&
	\bw^L_{i,j+\half} &= \bw_{i,j},& \bw^R_{i,j+\half} &= \bw_{i,j+\half}.
	\end{align*}
	For parabolic reconstruction we use a MUSCL scheme \cite{VanLeer1979}. It is given by
	\begin{align*}
	\bw^L_{i+\half,j} &= \bw_{i,j} + \frac 1 4 \left((1-\kappa)(\bw_{i,j}-\bw_{i-1,j})+(1+\kappa)(\bw_{i+1,j}-\bw_{i,j}\right),\\ 
	\bw^R_{i+\half,j} &= \bw_{i+1,j} - \frac 1 4 \left((1+\kappa)(\bw_{i+1,j}-\bw_{i,j})+(1-\kappa)(\bw_{i+2,j}-\bw_{i+1,j}\right),\\
	\bw^L_{i,j+\half} &= \bw_{i,j} + \frac 1 4 \left((1-\kappa)(\bw_{i,j}-\bw_{i,j-1})+(1+\kappa)(\bw_{i,j+1}-\bw_{i,j}\right),\\ 
	\bw^R_{i,j=\half} &= \bw_{i,j+1} - \frac 1 4 \left((1+\kappa)(\bw_{i,j+1}-\bw_{i,j})+(1-\kappa)(\bw_{i,j+2}-\bw_{i,j+1}\right).
	\end{align*}
	We obtain the reconstructed conservative variables as described in \cref{sec:wbScheme_fv}.
	
	When the constant reconstruction is used the resulting overall scheme is first order accurate. The parabolic MUSCL reconstruction combined with our well-balanced method leads to a second order accurate scheme due to the second order discretization of the gravitational source term.
	
	\allowdisplaybreaks[0]
	
	\subsection{Time-Stepping}
	
	For the numerical experiments in this paper, we use an explicit, three stage Runge--Kutta scheme, which has been introduced in \cite{Shu1988}. It has the \emph{total variation diminishing property} (\eg \cite{Toro2009}, chapter 13) and is third-order accurate in time.
	The scheme can be written as
	\begin{align}
	\vec q_{i,j}^{(1)} &= \vec q_{i,j}^n-\Delta t \vec R_{i,j}(\vec q^n),\\
	\vec q_{i,j}^{(2)} &= \frac{3}{4}\vec q_{i,j}^n +\frac{1}{4}\vec q_{i,j}^{(1)} -\frac{1}{4}\Delta t \vec R_{i,j}(\vec q^{(1)}),\\
	\vec q_{i,j}^{n+1} &= \frac{1}{3}\vec q_{i,j}^n +\frac{2}{3}\vec q_{i,j}^{(2)} -\frac{2}{3}\Delta t \vec R_{i,j}(\vec q^{(2)}),
	\label{eq:RK3}
	\end{align}
	where $\vec R_{i,j}$ is the spatial residual which contains the sum of the numerical fluxes and the discretized source term. The size of the time step $\Delta t$ is limited by the CFL condition and given via
	\[ \Delta t = C_{\text{CFL}} \cdot\min_{(i,j) \text{ with }(x_i,y_j)\in\Omega}\left(\frac{
		\min_{(k,l)\in\{(i+1,j),(i-1,j),(i,j+1),(i,j-1)\}}\|\x_{i,j}-\x_{k,l}\|_2
	}{c_{i,j}+|\vec v_{i,j}|}\right). \]
	In the numerical tests we use $C_\text{CFL}=0.9$. This choice leads to a CFL stable scheme \cite{Shu1988}.

	\subsection{Curvilinear Grids}
	\label{sec:slhgrids}
	
\begin{figure}
		\centering
		\subfloat [Cartesian grid
		\label{fig:cartesian}]
		{\includegraphics[width=0.25\textwidth]{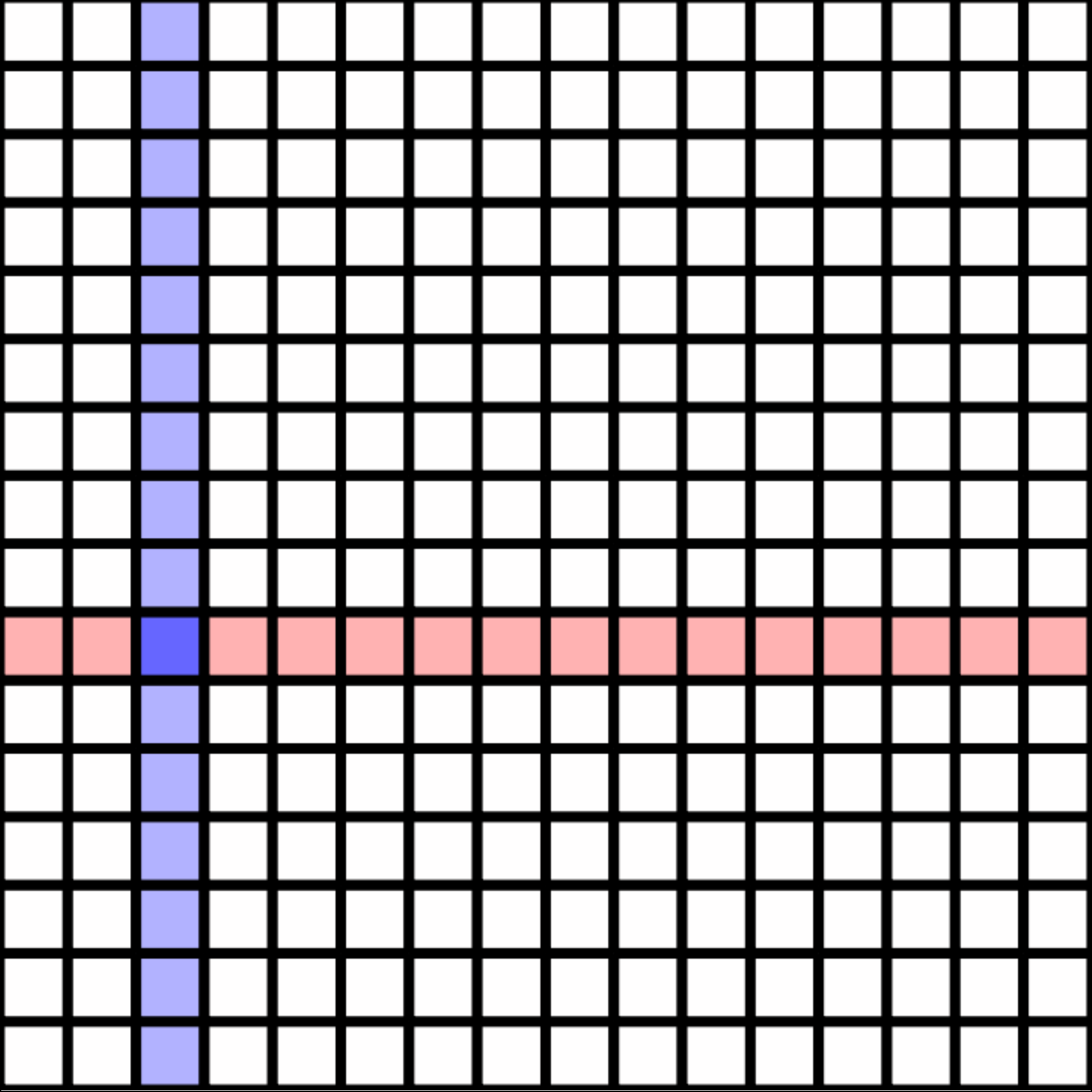}}
		\hfill
		\subfloat [Cubed sphere grid
		\label{fig:cubedsphere}]
		{\includegraphics[width=0.25\textwidth]{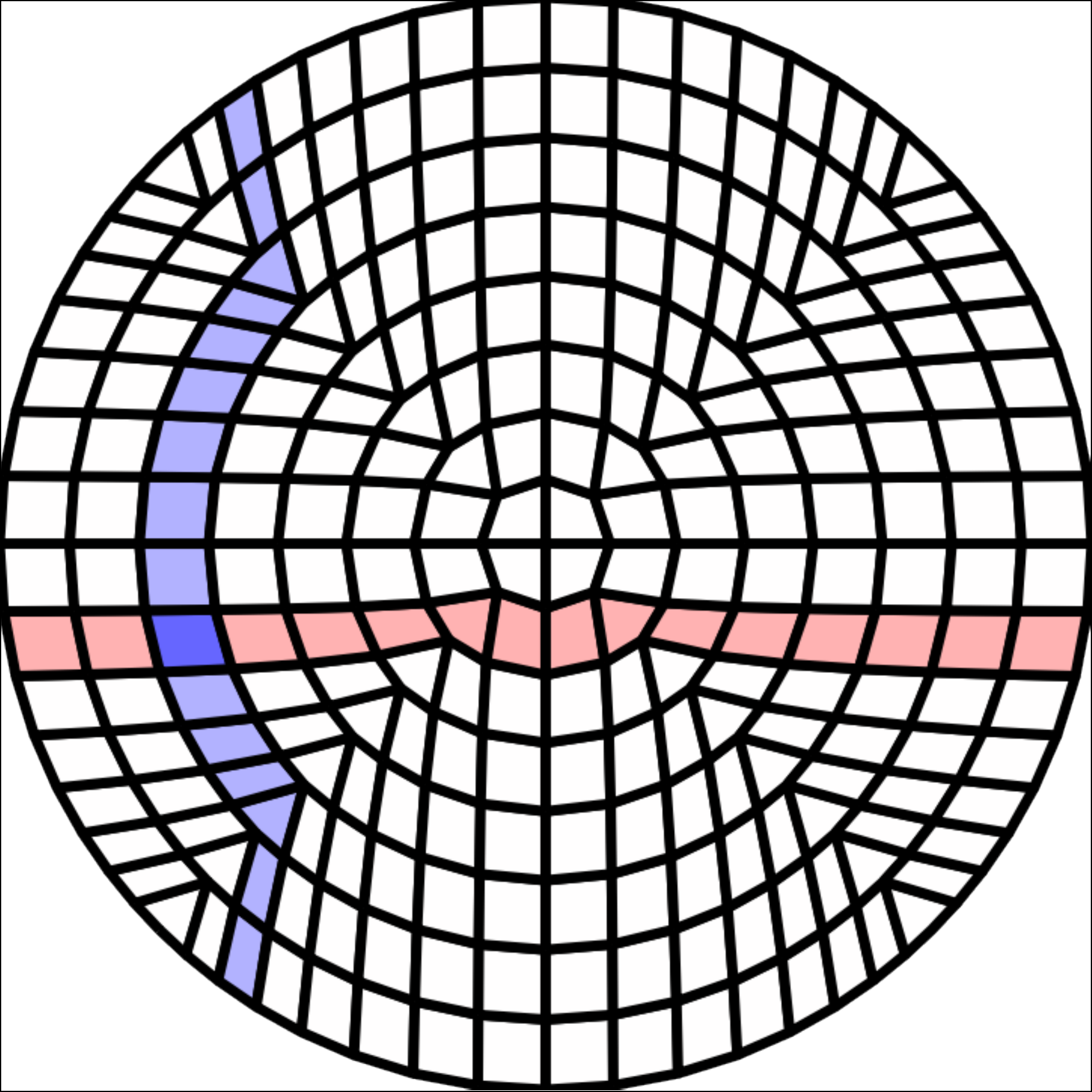}}
		\hfill
		\subfloat [Polar grid
		\label{fig:polar}]
		{\includegraphics[width=0.25\textwidth]{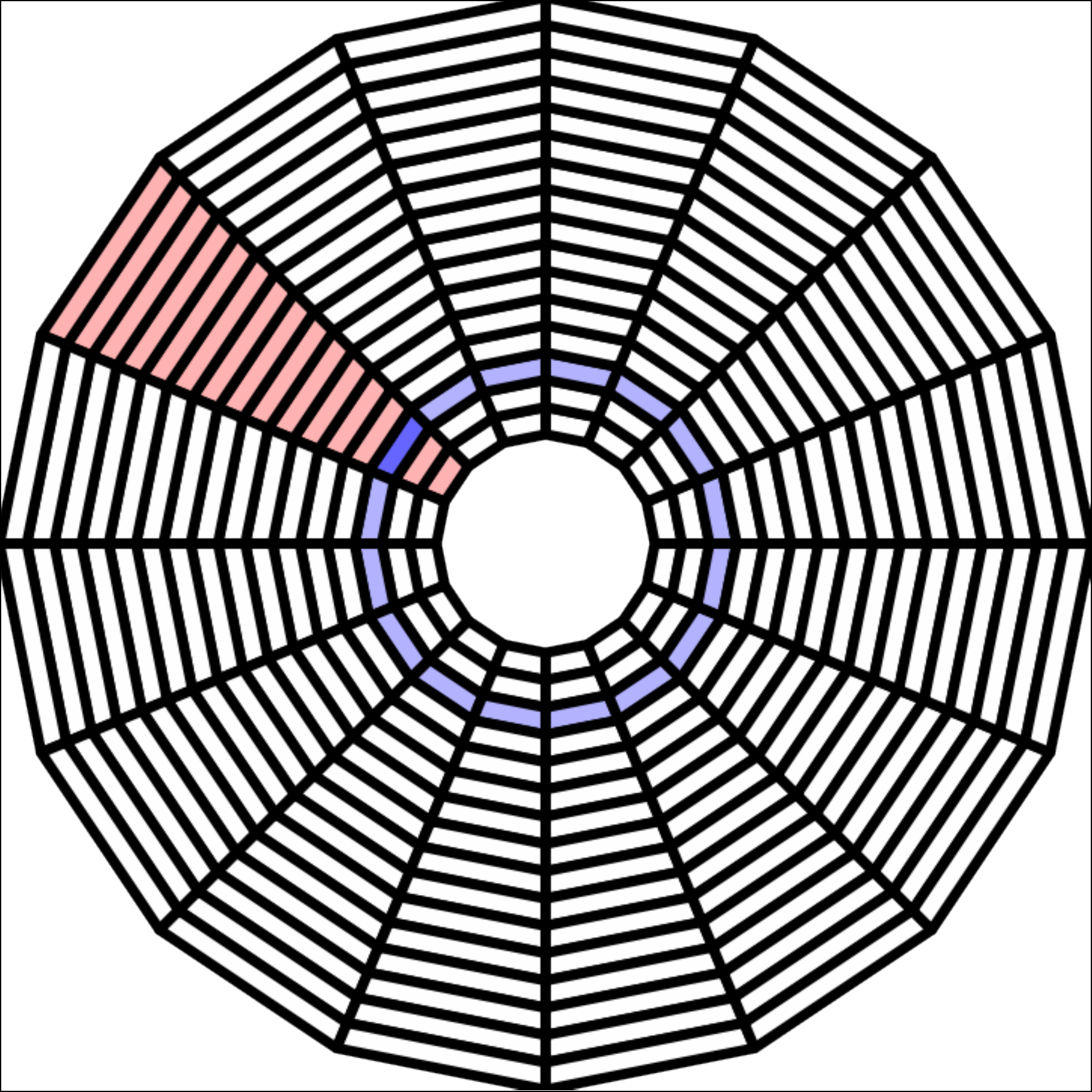}}
		\hfill
		\subfloat [Sinusoidal grid
		\label{fig:sinusoidal}]
		{\includegraphics[width=0.25\textwidth]{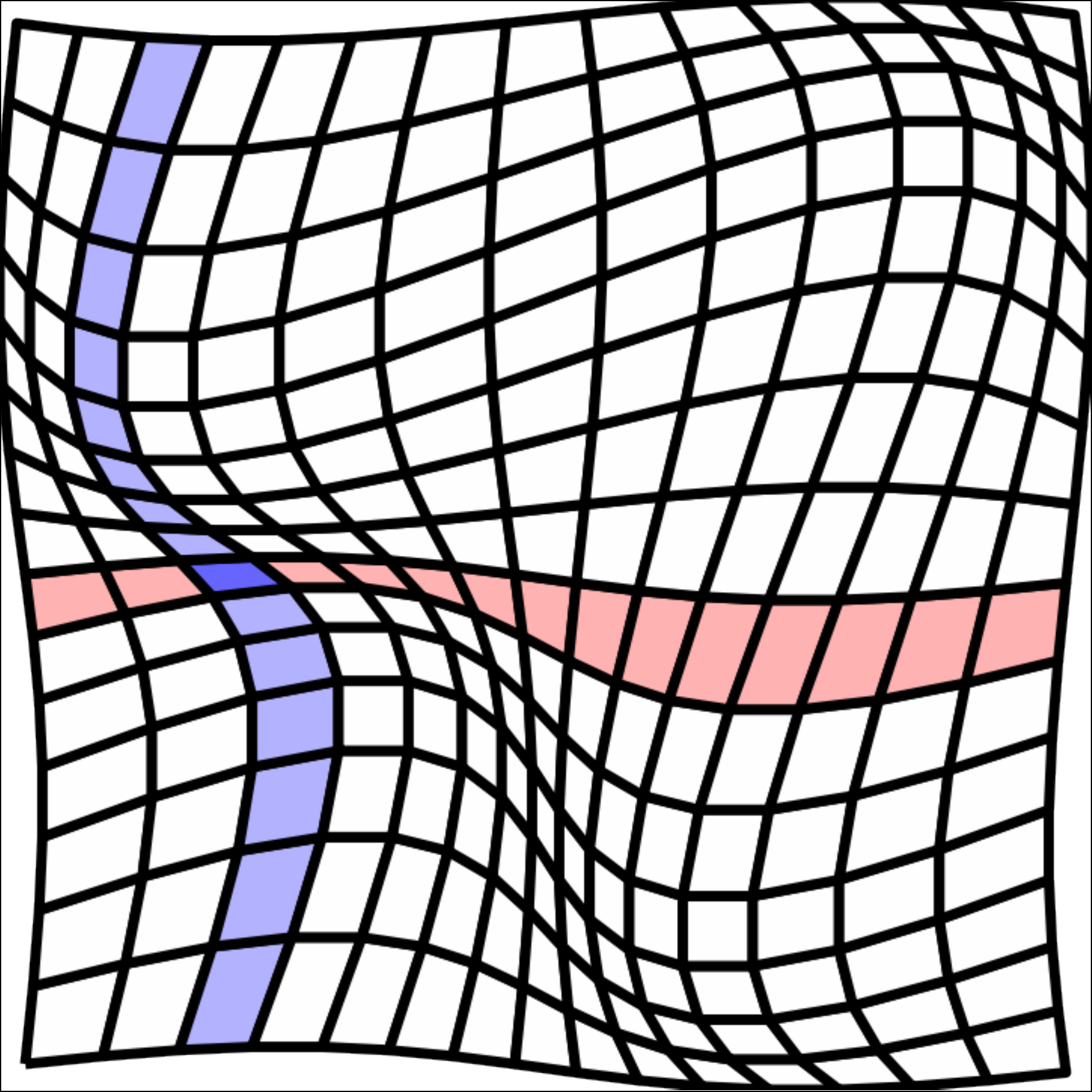}}
		\caption{\label{fig:grids}Some curvilinear grids which are used in this paper. Each of the grids is shown with $16\times16$ grid cells. The same row and line in the corresponding computational grid is colorized for every shown grid.}
	\end{figure}

	We will use curvilinear grids in many of the numerical tests. In practice, the transformed Euler equations \cref{eq:eul2d} are solved on a Cartesian grid on a reference domain.  Details on the construction of curvilinear grids are given in \cite{Kifonidis2012}, for example. 
	Three grids that are used in numerical tests in this paper are shown in \cref{fig:grids}. In \cref{fig:cartesian} we see the grid corresponding to the most trivial case of curvilinear coordinates, $\x(\boldsymbol\xi)=\boldsymbol\xi$. In \cref{fig:cubedsphere} the \emph{cubed sphere} grid is presented. It is the result of the attempt of Calhoun et.~al.~\cite{Calhoun2008} to create a grid that adapts to spherical symmetry, but is still basically Cartesian in the center. This way the central singularity of polar coordinates is avoided. A \emph{polar} grid, \ie a grid introduced by polar coordinates, can be seen in \cref{fig:polar}. It is discretely spherically symmetric. The center has to be avoided, since it is singular. The volume of the cells increases from smaller to higher radii. The \emph{sinusoidal} grid is shown in \cref{fig:sinusoidal}. This grid has been introduced by Colella et.~al.~\cite{Colella2011} for testing purposes. 
	
	\subsection{Non-Well-Balanced Scheme}
	
	Comparison with a non-well-balanced scheme is given in some parts of \cref{sec:numerical} to outline the importance of well-balancing. The difference to the well-balanced scheme presented in this section is in the reconstructed variables and the source term discretization. Instead of the hydrostatic variables $\vec w$ the primitive variables
	\[
	\vec w_\text{prim} := \left[\rho,u,v,p\right]^T
	\]
	are reconstructed. The source terms are computed by evaluating the gravitational acceleration at the cell-center exactly, so that the discretized source term has the form
	\[
	\vec S^\text{non-WB}_{i,j}:=\left[ 0, \,\rho_{i,j}\DD{\phi}{x}(x_i,y_j) , \,\rho_{i,j}\DD{\phi}{y}(x_i,y_j),\,0\right]^T.
	\]

	\section{Numerical Experiments}
	\label{sec:numerical}
	The methods discussed in \cref{sec:wbScheme} are applied in the numerical tests. For that we use the Seven-League Hydro Code (SLH) (see www.slh-code.org), which is an astrophysical code described in \cite{Miczek2013} and \cite{Edelmann2014}. All tests are performed in double precision.
	
	\subsection{Numerical Tests of the Stability}
	
	To demonstrate the well-balanced property of a scheme it is in principle sufficient to run a test with a hydrostatic setup for few time step. Yet, it is another question if the method is leading to a stable scheme or not. To address this question we run all tests for a longer time.
	For the tests we present in this section we check if the underlying problem is stable w.r.t.\ convection. This way we can distinguish numerical instabilities from underlying instabilities of the solution. In case the setup is stable w.r.t.\ convection we run the tests for $10\,\tBV$. The \brunt time is a natural time scale corresponding to the stability of the atmosphere.
	We concentrate on two dimensional problems since schemes are more likely to introduce instabilities in multidimensional applications due to the higher number of degrees of freedom. All tests in this subsection are conducted with the Roe numerical flux function combined with constant reconstruction and explicit RK3 time-stepping.
	The tests are conducted on the unit square $\Omega=[0,1]\times[0,1]$ using a Cartesian grid. Dirichlet boundary conditions are used. If some of the parameters are altered, they are given explicitly.
	
	\subsubsection{Isothermal Hydrostatic Solution}
	\label{sec:2d_isothermal}
	Consider the isothermal hydrostatic solution \cref{eq:isothermal}. To make the results comparable to \cite{Xing2013} and \cite{Chandrashekar2015} we choose a slightly modified problem given by
	\begin{equation}
	\bar\rho(x,y) = \rho_0\exp(-\rho_0\phi(x,y)/p_0),\quad 
	\bar p(x,y) = p_0\exp(-\rho_0\phi(x,y)/p_0)
	\label{eq:2d_isothermal}
	\end{equation}
	with $\phi(x,y) = x+y$, $\rho_0 = 1.21$, and $p_0=1$. We apply our \wb scheme combined with the Roe numerical flux function to evolve this initial condition to the time $t=10\,\tBV=117.55$. 
	
	The errors in $L^1$ norm with respect to the initial condition are shown in \subtabref{tab:2d_isothermal}. All errors are close to machine precision when the well-balanced scheme is applied. 
	When the non well-balanced scheme is applied, the errors are significant.
	\newcommand{\twodheader}{\head{cells}   &\head{scheme }&\head\Lorho&\head\Lou    &  \head\Lov  & \head\Lop        \\}
	\setcounter{subcount}{1}
	\begin{table}
		\centering
		\begin{tabularx}{1.0\textwidth}{ C{.13\textwidth} C{.12\textwidth} C{.15\textwidth} C{.15\textwidth} C{.15\textwidth} C{.15\textwidth}  C{.15\textwidth}}
			\toprule
			\twodheader
			\midrule
			\multicolumn{5}{l}{\sublabel{tab:2d_isothermal} \subheadstyle{Isothermal solution, ideal gas}} 														\\
			$50\times 50 $  & WB    & 7.9328e-15 & 3.0617e-16 & 4.3988e-15 & 2.1249e-15 \\
			& no WB & 1.3930e-02 & 6.6988e-15 & 8.6857e-03 & 3.9785e-03 \\
			$200\times 200$ & WB    & 5.3619e-14 & 5.5781e-16 & 2.8956e-14 & 1.1549e-14 \\
			& no WB & 3.4358e-03 & 7.0805e-15 & 2.1278e-03 & 9.0645e-04 \\
			\subheadseparator
			\multicolumn{5}{l}{\sublabel{tab:2d_isentropic} \subheadstyle{Isentropic solution, ideal gas}} 														\\
			$50\times 50 $  & WB    & 1.7171e-14 & 1.1442e-15 & 1.1455e-14 & 1.7990e-15 \\
			& no WB & 1.8785e-03 & 4.4098e-15 & 1.5233e-03 & 2.2894e-03 \\
			$200\times 200$ & WB    & 1.5001e-13 & 3.8162e-15 & 9.6811e-14 & 1.6175e-14 \\
			& no WB & 4.6335e-04 & 9.8775e-15 & 3.8098e-04 & 5.6391e-04 \\
			\subheadseparator
			\multicolumn{5}{l}{\sublabel{tab:2d_polytropic} \subheadstyle{Polytropic solution, $\nu=1.2<\gamma$, ideal gas}} 									\\
			$50\times 50 $  & WB    & 1.2129e-14 & 5.5438e-16 & 7.3212e-15 & 1.3942e-15 \\
			& no WB & 1.4198e-02 & 9.3538e-15 & 9.0315e-03 & 4.2003e-03 \\
			$200\times 200$ & WB    & 1.2081e-13 & 1.0600e-15 & 6.0344e-14 & 1.1861e-14 \\
			& no WB & 3.5177e-03 & 1.0028e-14 & 2.2226e-03 & 9.5553e-04 \\
			\subheadseparator
			\multicolumn{5}{l}{\sublabel{tab:2d_unstable_polytropic} \subheadstyle{Polytropic solution, $\nu=1.6>\gamma$, ideal gas}} 			\\
			$50\times 50 $  & WB    & 1.2100e-06 & 5.1134e-07 & 3.9891e-06 & 8.2899e-08 \\
			& no WB & 3.1469e-02 & 1.2112e-15 & 2.3118e-02 & 8.3768e-03 \\
			$200\times 200$ & WB    & 2.2515e-02 & 2.8865e-02 & 5.0755e-02 & 8.4021e-03 \\
			& no WB & 3.2213e-02 & 7.9883e-14 & 2.2901e-02 & 8.8411e-03 \\
			\subheadseparator
			\multicolumn{5}{l}{\sublabel{tab:2d_tanh} \subheadstyle{tanh-profile, ideal gas}} 																		\\
			$50\times 50 $  & WB    & 4.3523e-15 & 1.5694e-16 & 1.1054e-15 & 2.0373e-15 \\
			& no WB & 1.9732e-02 & 1.6283e-15 & 1.6307e-03 & 1.5402e-03 \\
			$200\times 200$ & WB    & 4.4862e-14 & 5.0936e-16 & 4.2709e-15 & 9.8190e-15 \\
			& no WB & 4.9237e-03 & 1.6620e-15 & 4.0395e-04 & 3.8237e-04 \\
			\subheadseparator
			\multicolumn{5}{l}{\sublabel{tab:2d_radiative_isothermal} \subheadstyle{Isothermal solution, ideal gas with radiation pressure}} 		 \\
			$50\times 50 $  & WB    & 2.0389e-14 & 4.0172e-16 & 6.8921e-15 & 1.1469e-14 \\
			& no WB & 1.3934e-02 & 1.3761e-14 & 5.3882e-03 & 4.7652e-03 \\
			$200\times 200$ & WB    & 1.5683e-13 & 7.6051e-16 & 2.8556e-14 & 5.5887e-14 \\
			& no WB & 3.4350e-03 & 1.4787e-14 & 1.3190e-03 & 9.4374e-04 \\
			\subheadseparator
			\multicolumn{7}{l}{\sublabel{tab:2d_radiative_polytropic} \subheadstyle{Polytropic solution, ideal gas with radiation pressure}}      \\
			$50\times 50 $  & WB    & 2.4695e-14 & 6.5524e-16 & 9.2055e-15 & 2.3608e-15 \\
			& no WB & 1.4167e-02 & 7.3164e-15 & 9.0282e-03 & 4.1308e-03 \\
			$200\times 200$ & WB    & 2.0169e-13 & 1.0721e-15 & 6.8032e-14 & 1.8171e-14 \\
			& no WB & 3.5143e-03 & 7.7572e-15 & 2.2251e-03 & 9.5087e-04 \\
			\subheadseparator
			\multicolumn{7}{l}{\sublabel{tab:2d_radiative_integrated} \subheadstyle{Numerically approximated solution, ideal gas with radiation pressure}} \\
			$50\times 50 $  & WB    & 6.0292e-14 & 1.1557e-14 & 1.1557e-14 & 9.6190e-15 \\
			& no WB & 1.0900e-02 & 5.9842e-03 & 5.9842e-03 & 4.2731e-03 \\
			$200\times 200$ & WB    & 3.3355e-13 & 4.8709e-14 & 4.8709e-14 & 5.1332e-14 \\
			& no WB & 3.0544e-03 & 1.5193e-03 & 1.5193e-03 & 1.2153e-03 \\
			\bottomrule
		\end{tabularx}
		\caption{$L^1$ errors for the tests in \crefrange{sec:2d_isothermal}{sec:rad_polytropic}. WB is short notation for well-balanced. The final time is $t=150$ for (b) and (d), $t=10\,\tBV$ for the other tests.}
		\label{tab:short_time}
	\end{table}
	
	\subsubsection{Isentropic Hydrostatic Solution}
	\label{sec:2d_isentropic}
	
	Consider the isentropic hydrostatic solution \cref{eq:polytropic} with $\nu=\gamma=1.4$. We use it as initial condition on the unit square with the gravitational potential $\phi(x,y) = x+y$. We run the simulation on different Cartesian grids to the final time $t=150$.  The results are shown in \subtabref{tab:2d_isentropic} using the $L^1$ norm. Like in \cref{sec:2d_isothermal} the well-balanced scheme keeps errors close to machine precision.

	\subsubsection{Polytropic Hydrostatic Solution ($\nu<\gamma$)}
	\label{sec:2d_polytropic}
	
\begin{figure}
		\centering
		\includegraphics[scale=1.]{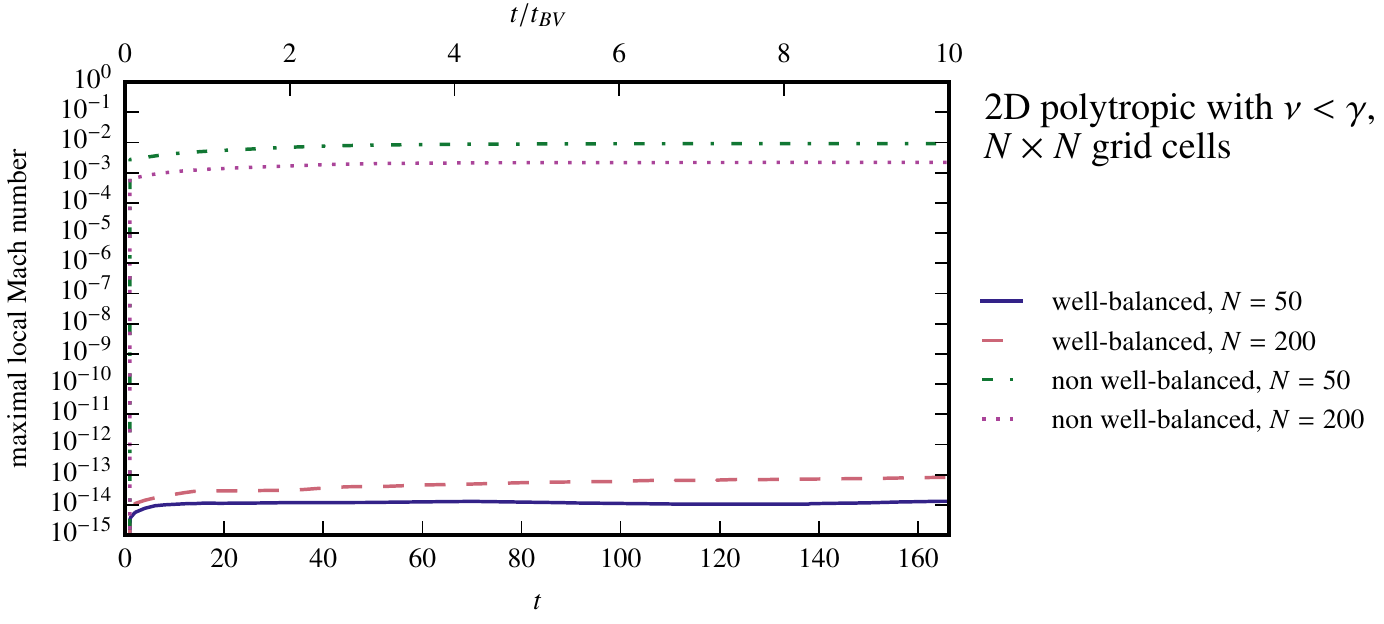}
		\caption{Maximal local Mach number (solid) and horizontal pressure fluctuations (dotted) for the two-dimensional polytropic equilibrium from \cref{sec:2d_polytropic} over time for $200\times200$ cells.}
		\label{fig:2d_polytropic}
	\end{figure}

	Consider the polytropic hydrostatic solution \cref{eq:polytropic} with an ideal gas EoS \cref{eq:ideal_eos} and $\nu=1.2$, \ie $\nu<\gamma$. We use it as initial condition on the unit square with the gravitational potential $\phi(x,y) = x+y$. We run the simulation on different Cartesian grids to the final time $t=10\,\tBV=166.24$.  The results are shown in \subtabref{tab:2d_polytropic} using the $L^1$ norm. Again, the well-balanced scheme keeps errors close to machine precision. The temporal evolution of the maximal local Mach number  can be seen in \cref{fig:2d_polytropic}

	\subsubsection{Polytropic Hydrostatic Solution ($\nu>\gamma$)}
	\label{sec:2d_unstable_polytropic}
	
\begin{figure}
		\centering
		\includegraphics[scale=1.]{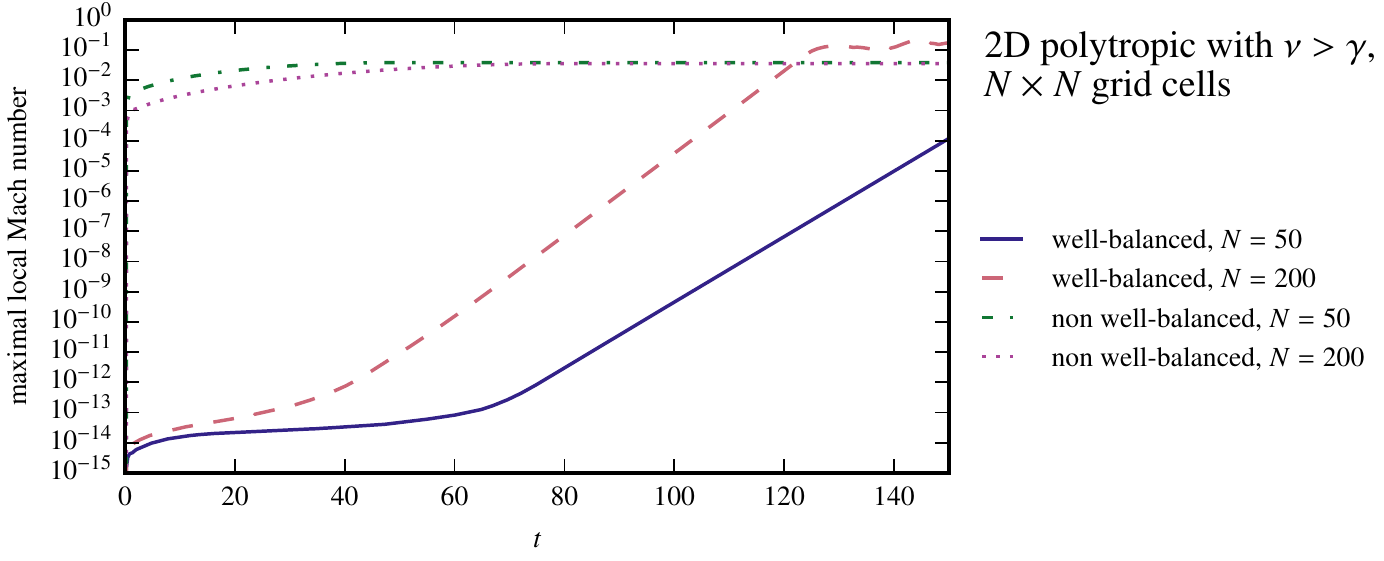}
		\caption{Maximal local Mach number (solid) and horizontal pressure fluctuations (dotted) for the two-dimensional solution which is unstable \wrt convection from \cref{sec:2d_unstable_polytropic} over time for $200\times200$ cells.}
		\label{fig:2d_unstable_polytropic}
	\end{figure}

	We redo the polytropic tests from \cref{sec:2d_polytropic} with $\nu=1.6$, \ie $\nu>\gamma$. The final time is $t=150$ since there is no \brunt time existent for $\nu>\gamma$. 
	The results are shown in \subtabref{tab:2d_unstable_polytropic} using the $L^1$ norm. The errors are significantly above machine precision, even if the well-balanced scheme is used. This is an expected result since the polytropic solution with $\nu>\gamma$ is not stable with respect to convection.
	The temporal evolution of the maximal local Mach number  can be seen in \cref{fig:2d_unstable_polytropic}
	
	\subsubsection{The tanh-Test}
	\label{sec:2d_tanh}
	
	Consider the hydrostatic solution \cref{eq:tanh_rhop}.
	In accordance to \cite{Edelmann2014} we choose $\Delta T= 0.1$, $\mu=0.02$, and the computational domain $x\in[-0.1,0.1]\times[-0.1,0.1]$.  For the \brunt time we obtain $\tBV=2.733$.
	We evolve these hydrostatic initial data to the final time $t=10\,\tBV=27.33$ on different Cartesian grids using our well-balanced scheme combined with the Roe numerical flux. The $L^1$ errors at final time are shown in \subtabref{tab:2d_tanh}. All errors are close to machine precision if the well-balanced scheme is used.

	\subsubsection{Ideal Gas with Radiation Pressure - Isothermal Solution}
	\label{sec:rad_isothermal}
	
	Let us consider the isothermal solution \cref{eq:radiative_isothermal} of \cref{eq:hystat} combined with the EoS \cref{eq:radiative_eos} with the diagonal gravitational potential $\phi(x,y)=x+y$. From \cref{sec:isothermal} we know that this setup is stable w.r.t\ convection. Using \cref{eq:brunt_rad_isothermal,eq:brunt_time} we compute the \brunt time $\tBV=24.657$. The $L^1$errors at $t=10\,\tBV=246.57$ can be seen in \subtabref{tab:2d_radiative_isothermal}. When the well-balanced scheme is applied, all errors are sufficiently close to machine precision.

	\subsubsection{Ideal Gas with Radiation Pressure - Polytropic Solution}
	\label{sec:rad_polytropic}
	We redo the tests from \cref{sec:2d_polytropic} with the EoS \cref{eq:radiative_eos} for an ideal gas with radiative pressure. As before, we ensure stability w.r.t.\ convection with a positive result and approximate the \brunt time $\tBV=12.885$ from the data using \cref{eq:brunt_ad,eq:ad_radiative,eq:brunt_time}. These initial data are evolved to $t=10\,\tBV=128.85$. The $L^1$ errors at final time are presented in \subtabref{tab:2d_radiative_polytropic}. For the \wb scheme, all errors are close to machine precision.
	
	\subsubsection{Ideal Gas with Radiation Pressure - Numerically Approximated Solution}
	\label{sec:rad_integrated}
	We run the test case with a the numerically approximated hydrostatic solution from \cref{sec:integratedProfile}. The data are computed using MATLAB. They are then given to our finite volume code in form of a table with discrete data points. The \brunt time $\tBV=25.708$ is computed from the data numerically using \cref{eq:brunt_ad,eq:ad_radiative}. We run the test to a final time of $t=10\,\tBV=257.08$. The $L^1$-errors at the final time are shown in \subtabref{tab:2d_radiative_integrated}. When the well-balanced scheme is applied, all errors are sufficiently close to machine precision.
	
	\subsection{Radial Test with different Grid Geometries}
	\label{sec:radialIsothermal}

	\def \myscale {1.}
	\begin{figure}[htpb]
		\centering
		\subfloat [Cartesian grid, \hfill\mbox{}\newline \mbox{}\hspace{.35cm} well-balanced\hfill\mbox{}
		\label{fig:radial_isothermal_0}]
		{\includegraphics[scale=\myscale]{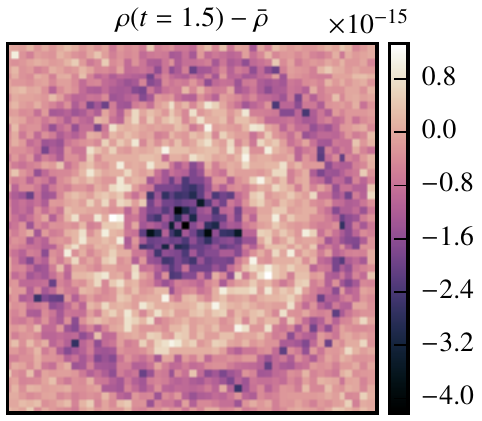}}
		\hfill
		\subfloat [cubed sphere grid, \hfill\mbox{}\newline\mbox{}\hspace{.35cm} well-balanced\hfill\mbox{}
		\label{fig:radial_isothermal_cs_0}]
		{\includegraphics[scale=\myscale]{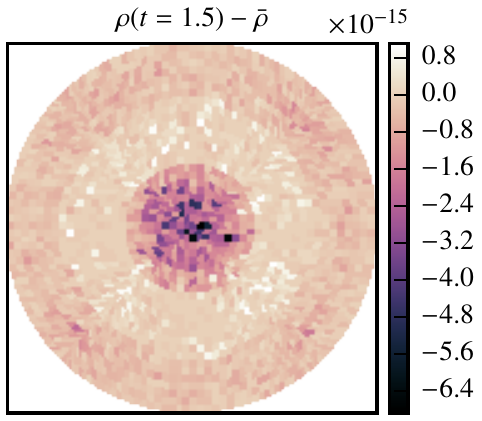}}
		\hfill
		\subfloat [polar grid,\hfill\mbox{} \newline\mbox{}\hspace{.35cm} well-balanced\hfill\mbox{}
		\label{fig:radial_isothermal_polar_0}]
		{\includegraphics[scale=\myscale]{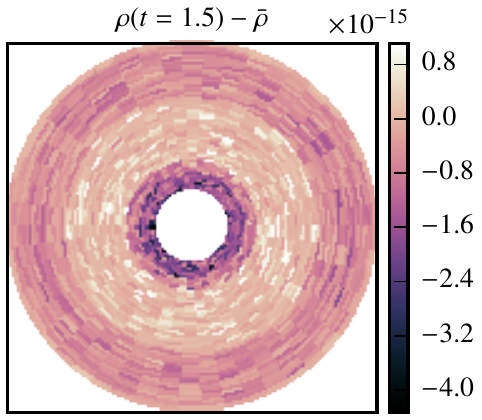}}
		
		\subfloat [Cartesian grid, \hfill\mbox{}\newline\mbox{}\hspace{.35cm} non-well-balanced\hfill\mbox{}
		\label{fig:radial_isothermal_2}]
		{\includegraphics[scale=\myscale]{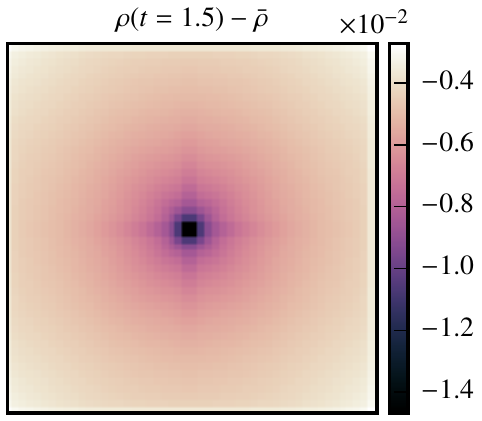}}
		\hfill
		\subfloat [cubed sphere grid, \hfill\mbox{}\newline\mbox{}\hspace{.35cm} non-well-balanced\hfill\mbox{}
		\label{fig:radial_isothermal_cs_2}]
		{\includegraphics[scale=\myscale]{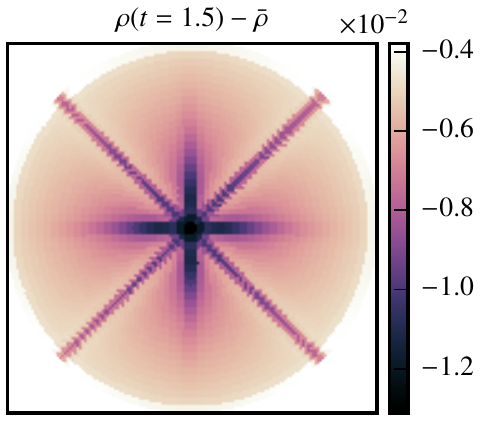}}
		\hfill
		\subfloat [polar grid, \hfill\mbox{}\newline\mbox{}\hspace{.35cm} non-well-balanced\hfill\mbox{}
		\label{fig:radial_isothermal_polar_2}]
		{\includegraphics[scale=\myscale]{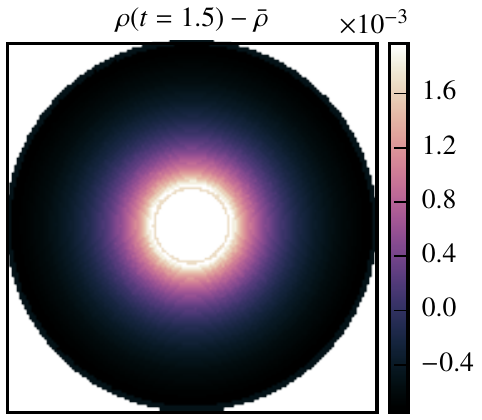}}
		\caption{%
			\label{fig:radial_isothermal}%
			Difference of the density from the initial condition at time $t=1.5$ for the isothermal test with radial potential from \cref{sec:radialIsothermal}.
			The resolution is $50\times50$ cells. The depicted domain is $\Omega=[-1,1]^2$. These figures are shown in the x-y plane.}
	\end{figure}
	
	We test the well-balanced property numerically for the isothermal solution \cref{eq:isothermal} for a radial gravitational potential $\phi(x,y)=\sqrt{x^2+y^2}$. The ideal gas EoS \cref{eq:ideal_eos} is used. The test is computed on a $50\times50$ Cartesian grid.
	The same configuration has been used for a well-balanced test in \cite{Chandrashekar2015}.
	As a second test with curvilinear grid, we use the cubed sphere grid suggested in \cite{Calhoun2008} with $50\times50$ cells. The structure of the grid is shown in \cref{fig:cubedsphere}.
	A polar grid as shown in \cref{fig:polar} is used as a third grid. For the tests with the polar grid we choose Dirichlet boundary conditions in the radial direction, periodic boundary conditions in the angular direction. The Cartesian and cubed sphere grids are combined with Dirichlet boundary conditions in both spatial directions.
	The domain is $\Omega=[-1,1]\times[-1,1]$ for the Cartesian grid. For the cubed sphere grid it reduces to the unit disk and for the polar grid we additionally subtract the disk around the origin with radius $0.2$. The test is conducted on all three grids with combinations of the well-balanced and the non-well-balanced scheme with the Roe numerical flux function and constant reconstruction. 
	
	The absolute density deviations from the initial data at time $t=1.5$ are visualized in \cref{fig:radial_isothermal}.
	In all cases the density errors of the well-balanced tests are close to machine precision while being of magnitude $10^{-2}$ to $10^{-3}$ in the non-well-balanced tests. In \cref{fig:radial_isothermal_cs_2} the the lack of rotation symmetry of the curvilinear grid in \cref{fig:cubedsphere} gets evident.
	
	\subsection{Evolution of a Small Perturbation}
	\label{sec:2d_perturbation}
	
	\def \perturbationscale {1.0}
	\begin{figure}
		{\bf Perturbation on an isothermal solution \cref{eq:2d_isothermal}}\\
		\subfloat[$\eta=10^{-1}$, \nwb\hfill\mbox{}
		\label{subfig:2d_iso_perturbation_2}]{%
			\includegraphics[scale=\perturbationscale]
			{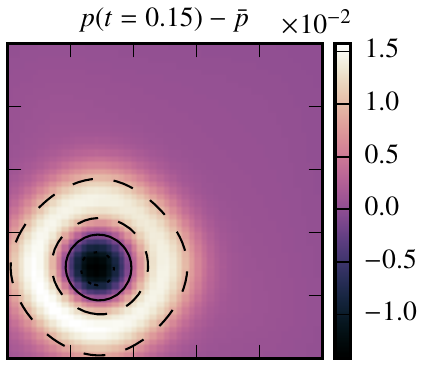}
		}\hfill
		\subfloat[$\eta=10^{-1}$, \wb\hfill\mbox{}
		\label{subfig:2d_iso_perturbation_1}]{%
			\includegraphics[scale=\perturbationscale]
			{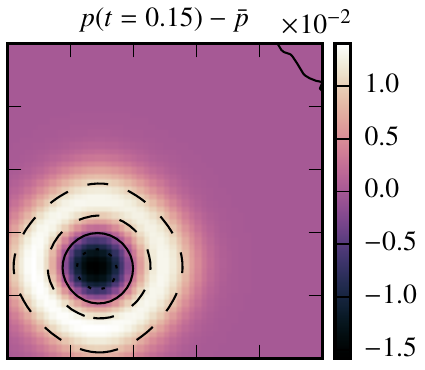}
		}\hfill
		\subfloat[$\eta=10^{-1}$, \wb,\hfill\mbox{} \newline\mbox{}\hspace{.35cm} sinusoidal grid
		\label{subfig:2d_iso_perturbation_5}]{%
			\includegraphics[scale=\perturbationscale]
			{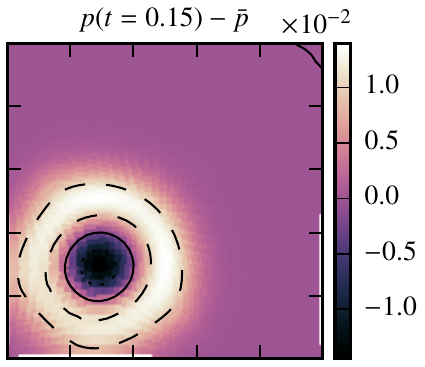}
		}
		
		\subfloat[$\eta=10^{-3}$, \nwb\hfill\mbox{}
		\label{subfig:2d_iso_perturbation_4}]{%
			\includegraphics[scale=\perturbationscale]
			{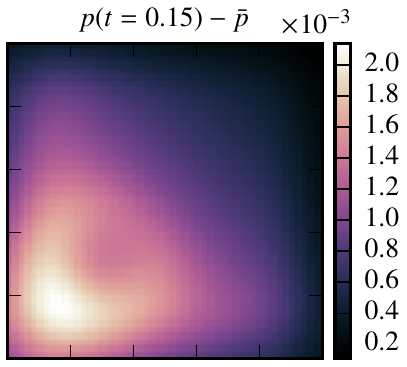}
		}\hfill
		\subfloat[$\eta=10^{-3}$, \wb\hfill\mbox{}
		\label{subfig:2d_iso_perturbation_3}]{%
			\includegraphics[scale=\perturbationscale]
			{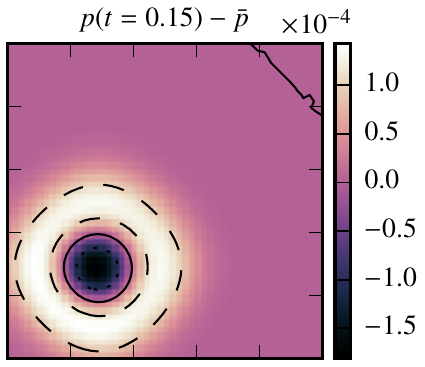}
		}\hfill
		\subfloat[$\eta=10^{-3}$, \wb, \hfill\mbox{}\newline\mbox{}\hspace{.35cm} sinusoidal grid
		\label{subfig:2d_iso_perturbation_6}]{%
			\includegraphics[scale=\perturbationscale]
			{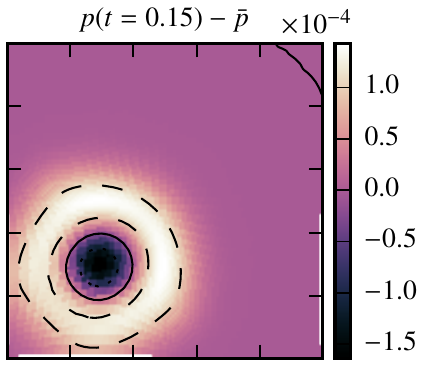}
		}\\[1.5em]
		
		{\bf Perturbation on a polytropic solution \cref{eq:polytropic} with $\nu=1.6>1.4=\gamma$}\\
		\subfloat[$\eta=10^{-1}$, \nwb\hfill\mbox{}
		\label{subfig:2d_poly_perturbation_2}]{%
			\includegraphics[scale=\perturbationscale]
			{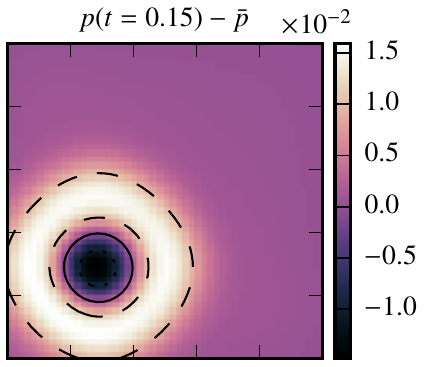}
		}\hfill
		\subfloat[$\eta=10^{-1}$, \wb\hfill\mbox{}
		\label{subfig:2d_poly_perturbation_1}]{%
			\includegraphics[scale=\perturbationscale]
			{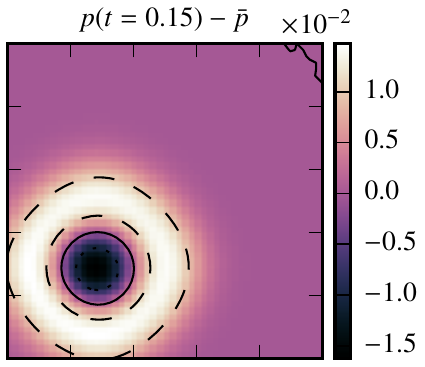}
		}\hfill
		\subfloat[$\eta=10^{-10}$, \wb\hfill\mbox{}
		\label{subfig:2d_poly_perturbation_5}]{%
			\includegraphics[scale=\perturbationscale]
			{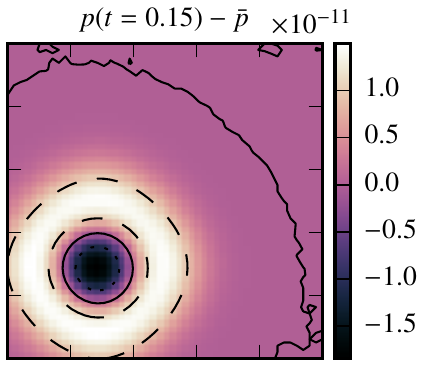}
		}

		\subfloat[$\eta=10^{-3}$, \nwb\hfill\mbox{}
		\label{subfig:2d_poly_perturbation_4}]{%
			\includegraphics[scale=\perturbationscale]
			{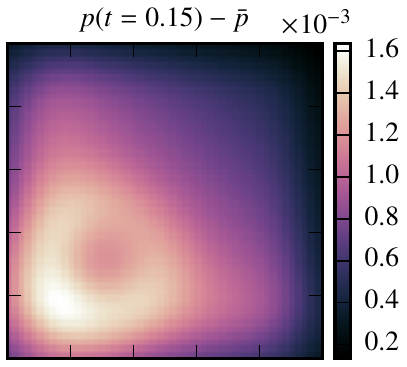}
		}\hfill
		\subfloat[$\eta=10^{-3}$, \wb\hfill\mbox{}
		\label{subfig:2d_poly_perturbation_3}]{%
			\includegraphics[scale=\perturbationscale]
			{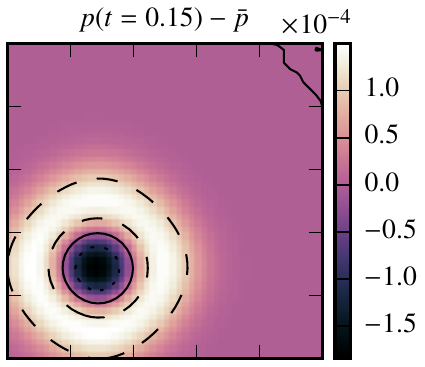}
		}\hfill
		\subfloat[$\eta=10^{-12}$, \wb\hfill\mbox{}
		\label{subfig:2d_poly_perturbation_6}]{%
			\includegraphics[scale=\perturbationscale]
			{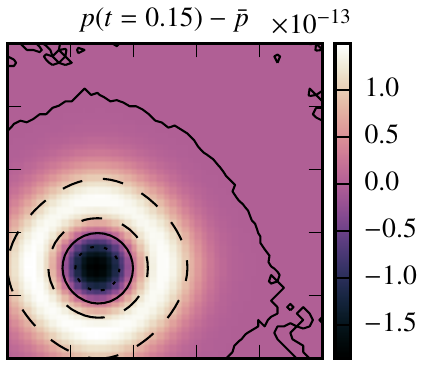}
		}
		\caption{Evolved pressure perturbation (\cref{eq:2d_perturbation}) on an isothermal and a polytropic hydrostatic solution at time $t=0.15$.
			The $x$-coordinate increases to the right, the $y$-coordinate to the top.
			Lines: dotted at $-0.1\eta$, solid at $0.0$, dashed at $0.1\eta$.}
		\label{fig:2d_perturbation}
	\end{figure}
	
	\begin{figure}
		\centering
		{\bf Perturbation on the numerically approximated hydrostatic solution\\ for ideal gas with radiation pressure from \cref{sec:integratedProfile}}\\
		\subfloat[$\eta=10^{-1}$, \nwb\hfill\mbox{}
		\label{subfig:2d_rad_perturbation_2}]{%
			\includegraphics[scale=\perturbationscale]
			{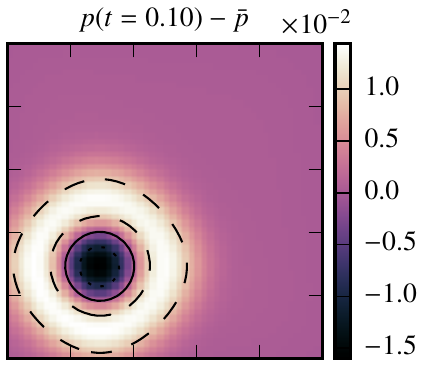}
		}\hfill
		\subfloat[$\eta=10^{-1}$, \wb\hfill\mbox{}
		\label{subfig:2d_rad_perturbation_1}]{%
			\includegraphics[scale=\perturbationscale]
			{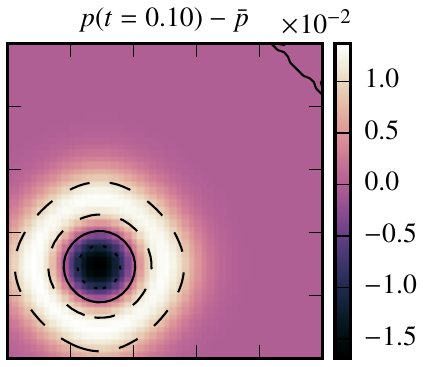}
		}\hfill
		\subfloat[$\eta=10^{-1}$, \nwb, \hfill\mbox{}\newline\mbox{}\hspace{.35cm} $1000\times1000$ grid
		\label{subfig:2d_rad_perturbation_5}]{%
			\includegraphics[scale=\perturbationscale]
			{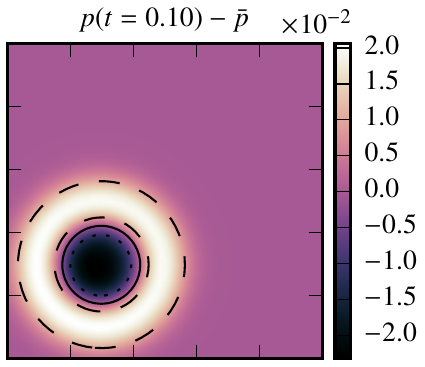}
		}
		
		\subfloat[$\eta=10^{-3}$, \nwb\hfill\mbox{}
		\label{subfig:2d_rad_perturbation_4}]{%
			\includegraphics[scale=\perturbationscale]
			{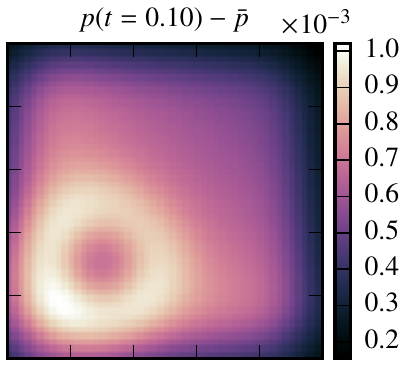}
		}\hfill
		\subfloat[$\eta=10^{-3}$, \wb\hfill\mbox{}
		\label{subfig:2d_rad_perturbation_3}]{%
			\includegraphics[scale=\perturbationscale]
			{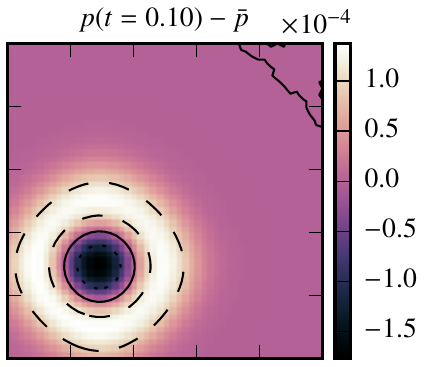}
		}\hfill
		\subfloat[$\eta=10^{-3}$, \nwb,\hfill\mbox{} \newline\mbox{}\hspace{.35cm} $1000\times1000$ grid
		\label{subfig:2d_rad_perturbation_6}]{%
			\includegraphics[scale=\perturbationscale]
			{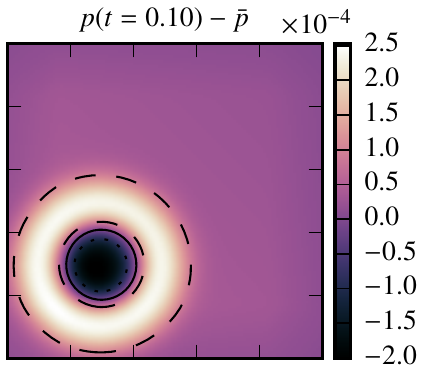}
		}
		\caption{Evolved pressure perturbation (\cref{eq:2d_perturbation}) on the numerically approximated hydrostatic solution  (\cref{sec:integratedProfile}) for ideal gas with radiation pressure at time $t=0.1$.
			The $x$-coordinate increases to the right, the $y$-coordinate to the top.
			Lines: dotted at $-0.1\eta$, solid at $0.0$, dashed at $0.1\eta$.}
		\label{fig:2d_perturbation_rad}
	\end{figure}
	\begin{figure}
		\centering
		\includegraphics[scale=\perturbationscale]{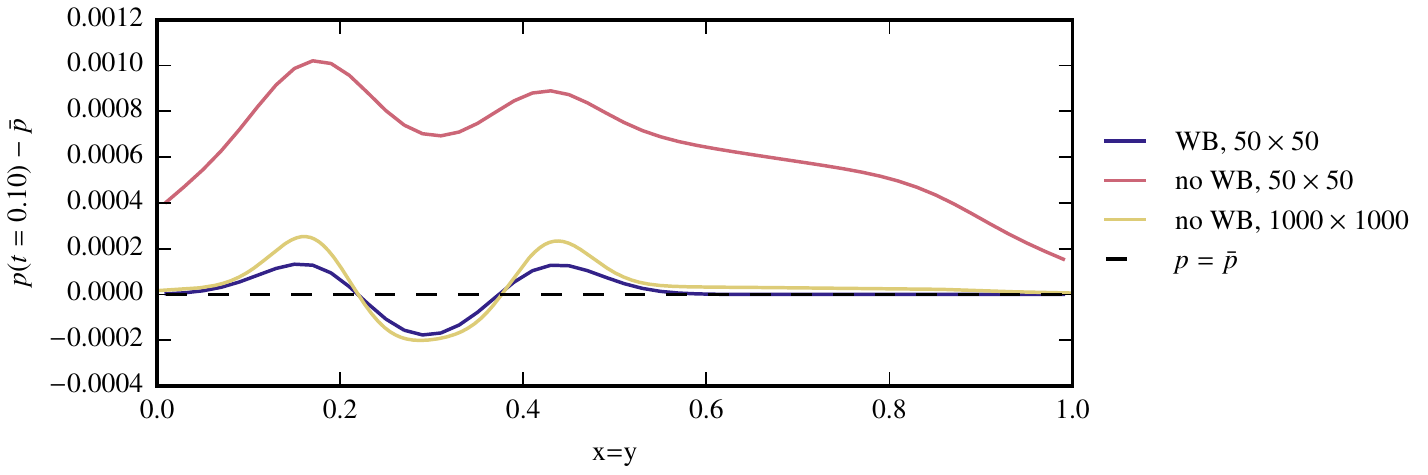}
		\caption{Diagonal cut along the line $x=y$ through the domains corresponding to $\eta=0.001$ in \cref{fig:2d_perturbation_rad}.}
		\label{fig:diagonalcut}
	\end{figure}
	In this test, we study the evolution of a small perturbation added to the hydrostatic solution. We test this for the isothermal solution \cref{eq:2d_isothermal}, the polytropic solution \cref{eq:polytropic} with an ideal gas EoS \cref{eq:ideal_eos} and $\nu=1.2$, \ie $\nu<\gamma$, and on the numerically approximated hydrostatic solution for ideal gas with radiation pressure (\cref{sec:integratedProfile}). In correspondence with \cite{Chandrashekar2015} and \cite{Xing2013}, we choose the same parameters as in the \cref{sec:2d_isothermal,sec:2d_polytropic} respectively and the initial pressure is
	\begin{equation}
	p(x,y,0) = \bar p(x,y) + \eta\exp(-100\rho_0((x-0.3)^2+(y-0.3)^2)/p_0).
	\label{eq:2d_perturbation}
	\end{equation}
	The initial density is $\rho(\cdot,\cdot,0)=\bar\rho$. The grid resolution $50\times50$. For the numerically approximated hydrostatic solution a $1000\times1000$ grid is used to compute reference solutions. The final time is $t=0.15$ for the tests on isothermal and polytropic backgrounds. It is $t=0.10$ for the tests with the numerically approximated hydrostatic solution as background.
	
	The pressure perturbations at final time can be seen in \cref{fig:2d_perturbation,fig:2d_perturbation_rad}. They are similar for all three equilibria: The large perturbation with $\eta=0.1$ is well-resolved for both, the \wb and the \nwb scheme. When the perturbation is decreased to $\eta=0.001$, the \nwb scheme is not able to resolve it well anymore, since the discretization errors start to dominate after some time. The \wb scheme shows no problems for the smaller perturbation.
	The isothermal test case has also been conducted on a sinusoidal grid. This grid is introduced in \cref{sec:slhgrids}. 
	The result of the test can be seen in \cref{subfig:2d_iso_perturbation_5,subfig:2d_iso_perturbation_6}. We see that the usage of the curvilinear mesh does not introduce significant errors in this test.
	In Figures (\ref{subfig:2d_poly_perturbation_5}) and (\ref{subfig:2d_poly_perturbation_6}) we can see that even a small perturbation of $\eta=10^{-10}$ or less leads to a well-resolved result, if the \wb scheme is used.
	
	In \cref{fig:diagonalcut} a cut through the domain along the diagonal $x=y$ is shown for the test with the numerically approximated hydrostatic solution. The pressure perturbation from the hydrostatic solution at the final time $t=0.1$ is shown along this line. The solution obtained with the well-balanced method on the coarse grid and the highly resolved solution with the standard scheme are close to each other. On the fine grid the diffusion is lower and the amplitude of the remaining perturbation is thus higher. In the test on the coarse grid without the well-balanced scheme the perturbation due to discretization errors dominates.

	\subsection{Radial Rayleigh-Taylor Instability}
	
	\begin{figure}[tb]
		\centering
		{\begin{center}\bf Cartesian grid\end{center}}
		\includegraphics[scale=1.]{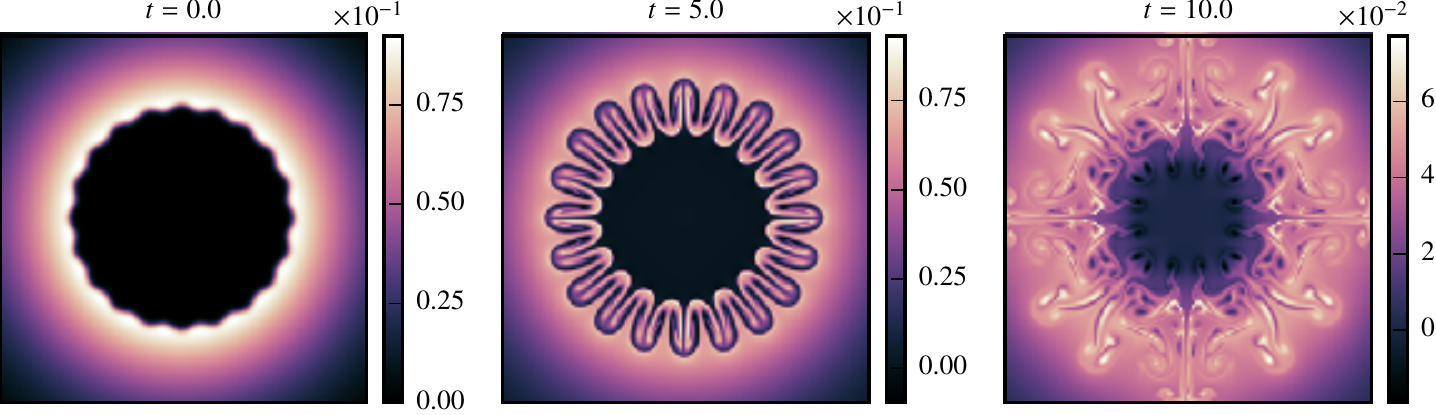}
		{\begin{center}\bf cubed sphere grid\end{center}}
		\includegraphics[scale=1.]{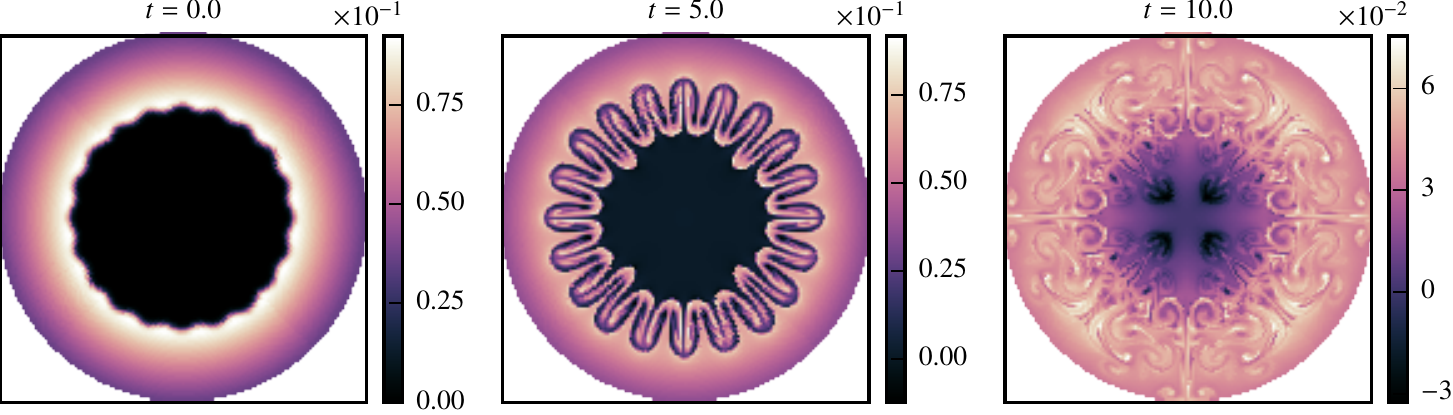}
		{\begin{center}\bf polar grid\end{center}}
		\includegraphics[scale=1.]{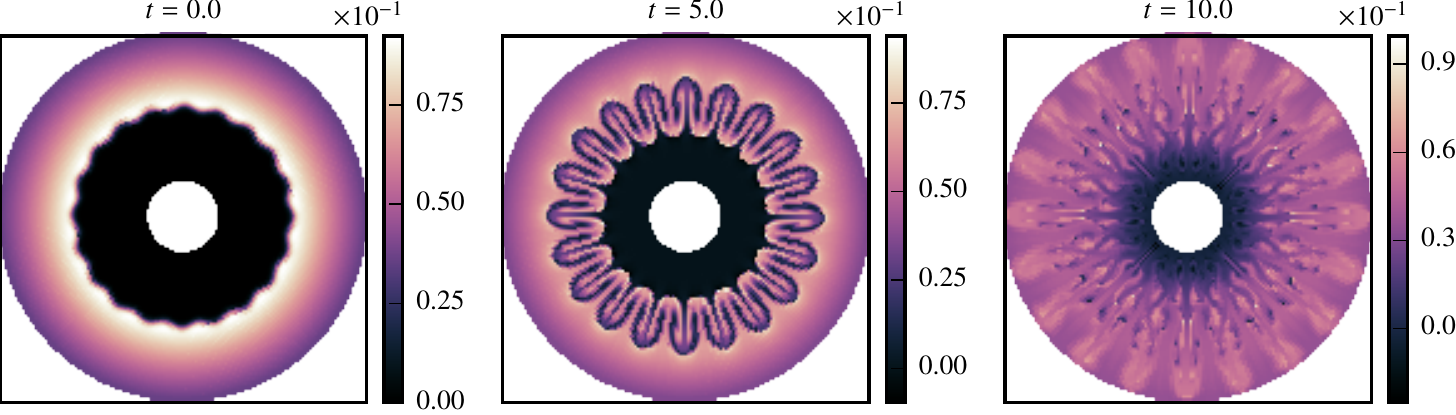}
		\caption{Temporal evolution of the density field for the radial Rayleigh--Taylor test. Our well-balanced scheme is applied in combination with the Roe numerical flux function, linear reconstruction, explicit RK3 time stepping, and different grid geometries. The grid size is $256\times256$ cells on the domain $[-1,1]^2$. The figures are shown in the x-y plane.
			The quantity shown is $\rho-\rho_{\text{internal}}$.}
		\label{fig:rrt}
	\end{figure}
	
	The next test case is given by a piecewise isothermal hydrostatic solution. The gravitational potential is given by $\phi(x,y)=r:=\sqrt{x^2+y^2}$, just as in \cref{sec:radialIsothermal}. The domain is $\Omega=[0,1]\times[0,1]$ and the initial condition is given by
	\begin{align}
	p&=\left\{
	\begin{array}{ll}
	\exp(-r) & \text{if } r\leq r_0\\
	\exp\left(\frac{-r+r_0(1-a)}{a}\right) &\text{if } r>r_0
	\end{array}
	\right. ,\\
	\rho&= \left\{
	\begin{array}{ll}
	\exp(-r) & \text{if } r\leq r_i\\
	\frac{1}{a}\exp\left(\frac{-r+r_0(1-a)}{a}\right) &\text{if } r>r_i
	\end{array}
	\right. ,
	\label{eq:rrt}
	\end{align}
	with
	\[
	r_i := r_0\left( 1+\eta\cos\left( k\theta \right) \right), \qquad
	a   := \frac{\exp(-r_0)}{\exp(-r_0)+\Delta\rho}.
	\]
	While the pressure is continuous on the whole domain, the density jumps at $r=r_i$ by approximately $\Delta\rho>0$.
	In accordance to \cite{Leveque1999}, \cite{Chandrashekar2015}, and \cite{Chandrashekar2017}, we choose the parameters $\Delta\rho=0.1$, $\eta=0.02$, and $k=20$, and the grid consists of $256\times256$ cells.
	As the solution evolves with time, mixing at the boundary between the equilibria is expected.
	
	In \cref{fig:rrt} results of the test are presented. Our well-balanced scheme is applied in combination with the Roe numerical flux function, MUSCL reconstruction, and explicit RK3 time stepping. The shown quantity is $\rho-\rho_\text{internal}$, where $\rho_{\text{internal}}(x,y):=\exp(-r(x,y))$. The functions $\alpha$ and $\beta$ are chosen to balance the internal hydrostatic solution on the whole grid. The test is run to a maximal time of $t=10.0$. Three different grid geometries have been used: a Cartesian grid, a cubed sphere grid, and a polar grid (see \cref{sec:slhgrids}), each with a solution of $256\times 256$ grid cells. At time t = 10 we show the mixing layer for all three grids.

	\subsection{Hot Rising Bubble}
	\label{sec:hotbubble}
	
	In this section we use a test case that has been used for example in \cite{Mendez-Nunez1994} and \cite{Rieper2013}. The problem consists of an isentropic atmosphere with a spot (bubble) of increased temperature. The initial data are given by
	\begin{align}
	\rho_\text{init}(x,y)&=\hat\rho\left( 1-\frac{\hat\rho(\gamma-1)}{\hat p \gamma}\phi(x,y) \right)^{\oo{\gamma-1}},\\
	p_\text{init}(x,y)&=\hat p^{1-\gamma}\left( R \theta(x,y)\rho_\text{init}(x,y) \right)^\gamma,\\
	u_\text{init}(x,y)&=0,\\
	v_\text{init}(x,y)&=0
	\label{eq:hotbubble}
	\end{align}
	with
	\begin{align}
	\phi(x,y)&=g y,\\
	r(x,y) &=\sqrt{(x-x_0)^{2}+(y-y_0)^2},\\
	\theta(x,y) &=\left\{
	\begin{array}{ll}
	\hat\theta\left(1+\Delta\theta \cos^2\left( \frac{\pi r(x,y)}{2 r_0} \right)\right) & \text{ if } r(x,y)\leq r_0\\
	\hat\theta & \text{ if } r(x,y)>r_0
	\end{array}
	\right.
	\label{eq:theta}
	\end{align}
	on the domain $\Omega=[0,10^6]\times[0,1.5\cdot10^6]$.
	For the constants we choose
	\newcommand{\mydist}{\quad\;\;\;}
	\begin{gather}
	R=8.314472\cdot 10^7,\mydist\hat\theta=300,\mydist\Delta\theta=0.022,\mydist \hat p=10^6,\mydist \hat\rho=\frac{\hat p}{\hat\theta R},
	\label{eq:bubbleconsts1}\\
	g=981,\mydist x_0=5\cdot10^5,\mydist y_0=2.75\cdot 10^5,\mydist\text{and}\mydist r_0=2.5\cdot 10^5.\label{eq:bubbleconsts2}
	\end{gather}
	We use the ideal gas EoS with the value for the gas constant $R$ given in \cref{eq:bubbleconsts1}. 
	The above values are typical for a terrestrial atmosphere in CGS units (\eg \cite{Kraemer2013}). 
	\begin{figure}
		\centering
		\includegraphics[scale=1]{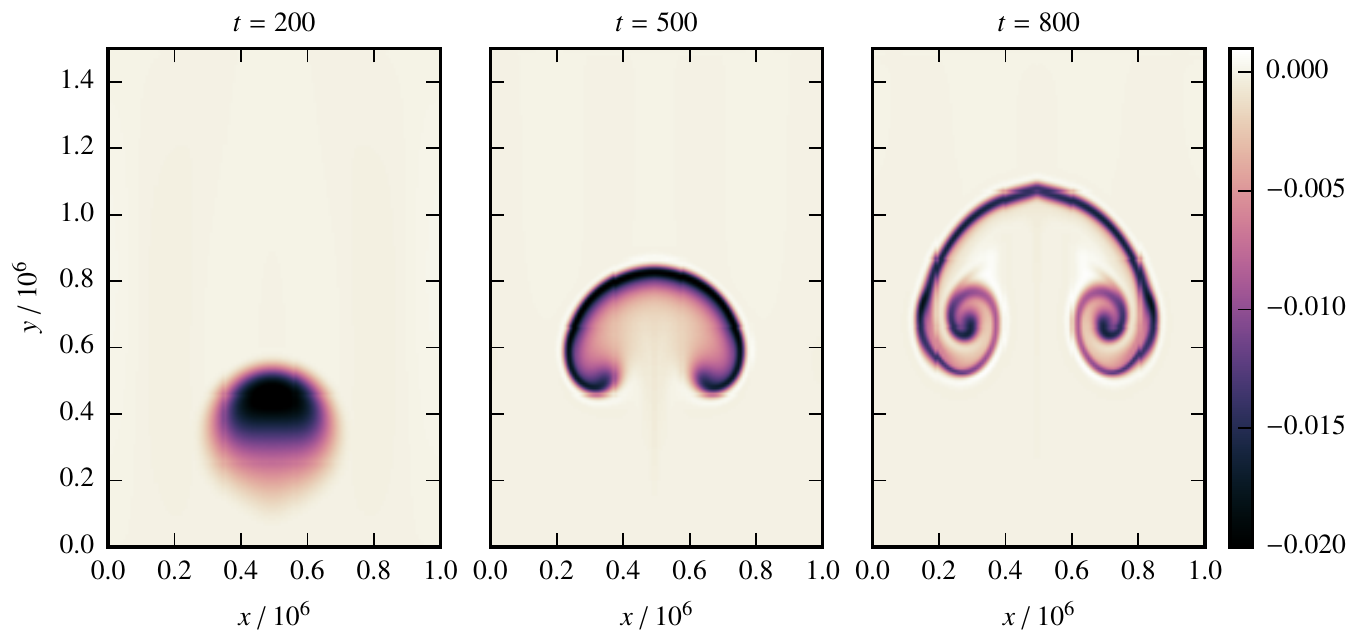}
		\caption{Relative density deviation for the hot rising bubble test in \cref{sec:hotbubble}.
			Our well-balanced scheme is applied in combination with the Roe numerical flux function, linear reconstruction, and explicit RK3 time stepping. The grid size is $128\times 192$ cells
		}
		\label{fig:hotbubble}
	\end{figure}
	
	Tests are conducted on a Cartesian $128\times192$ grid. Our well-balanced scheme is applied in combination with the Roe numerical flux function, MUSCL reconstruction, and explicit RK3 time stepping. We choose the functions $\alpha$ and $\beta$ in such a way that the isentropic background --\ie \crefrange{eq:hotbubble}{eq:bubbleconsts2} with $\Delta\theta=0$-- is balanced.
	Results of the tests are presented in \cref{fig:hotbubble}.

	\subsection{Order of Accuracy Study}
	\label{sec:accuracy}
	
	\begin{table}
		\centering
		\begin{tabular}{ccccccccc}
			\multicolumn{9}{l}{\bf{Cartesian grid:}}\\
			\toprule
			N & $\rho$: Error & Rate & $\rho u$: Error & Rate& $\rho v$: Error & Rate& $\rho E$: Error & Rate \\
			\midrule
			16 & 0.0006991 & -- & 0.00071139 &  -- & 0.00071139 &  -- & 0.001336 &  -- \\
			32 & 0.00013033 & 2.42 & 0.00013137 & 2.44 & 0.00013137 & 2.44 & 0.000237 & 2.49 \\
			64 & 2.6642e-05 & 2.29 & 2.7032e-05 & 2.28 & 2.7032e-05 & 2.28 & 4.6686e-05 & 2.34 \\
			128 & 5.9188e-06 & 2.17 & 6.0525e-06 & 2.16 & 6.0525e-06 & 2.16 & 1.011e-05 & 2.21 \\
			256 & 1.3883e-06 & 2.09 & 1.4286e-06 & 2.08 & 1.4286e-06 & 2.08 & 2.3413e-06 & 2.11 \\
			512 & 3.3579e-07 & 2.05 & 3.4684e-07 & 2.04 & 3.4684e-07 & 2.04 & 5.627e-07 & 2.06 \\
			1024 & 8.2542e-08 & 2.02 & 8.5441e-08 & 2.02 & 8.5441e-08 & 2.02 & 1.379e-07 & 2.03 \\
			\bottomrule
			\\
			\multicolumn{9}{l}{\bf{sinusoidal grid:}}\\
			\toprule
			N & $\rho$: Error & Rate & $\rho u$: Error & Rate& $\rho v$: Error & Rate& $\rho E$: Error & Rate \\
			\midrule
			16 & 0.0029866 & -- & 0.0058612 &  -- & 0.0058612 &  -- &  0.0636 & -- \\
			32 & 0.00067682 & 2.14 & 0.0015231 & 1.94 & 0.0015231 & 1.94 & 0.015522 & 2.03 \\
			64 & 0.00015966 & 2.08 & 0.0003877 & 1.97 & 0.0003877 & 1.97 & 0.0037432 & 2.05 \\
			128 & 3.875e-05 & 2.04 & 9.8001e-05 & 1.98 & 9.8001e-05 & 1.98 & 0.00091998 & 2.02 \\
			256 & 9.5467e-06 & 2.02 & 2.4647e-05 & 1.99 & 2.4647e-05 & 1.99 & 0.00022799 & 2.01 \\
			512 & 2.3709e-06 & 2.01 & 6.1805e-06 & 2.00 & 6.1805e-06 & 2.00 & 5.6744e-05 & 2.01 \\
			1024 & 5.9083e-07 & 2.00 & 1.5475e-06 & 2.00 & 1.5475e-06 & 2.00 & 1.4154e-05 & 2.00 \\
			\bottomrule
			\\
			\multicolumn{9}{l}{\bf{polar grid:}}\\
			\toprule
			N & $\rho$: Error & Rate & $\rho u$: Error & Rate& $\rho v$: Error & Rate& $\rho E$: Error & Rate \\
			\midrule
			16 & 0.015013 & -- & 0.01851 & -- & 0.01851 & -- & 0.057656 & -- \\
			32 & 0.0050899 & 1.56 & 0.006662 & 1.47 & 0.006662 & 1.47 & 0.016689 & 1.79 \\
			64 & 0.00097143 & 2.39 & 0.0013992 & 2.25 & 0.0013992 & 2.25 & 0.0045171 & 1.89 \\
			128 & 0.0001856 & 2.39 & 0.00030959 & 2.18 & 0.00030959 & 2.18 & 0.0011836 & 1.93 \\
			256 & 3.8076e-05 & 2.29 & 7.232e-05 & 2.10 & 7.232e-05 & 2.10 & 0.00030843 & 1.94 \\
			512 & 8.4367e-06 & 2.17 & 1.7567e-05 & 2.04 & 1.7567e-05 & 2.04 & 7.9112e-05 & 1.96 \\
			1024 & 1.9729e-06 & 2.10 & 4.3412e-06 & 2.02 & 4.3412e-06 & 2.02 & 2.0045e-05 & 1.98 \\
			\bottomrule
			\\
			\multicolumn{9}{l}{\bf{cubed sphere grid:}}\\
			\toprule
			N & $\rho$: Error & Rate & $\rho u$: Error & Rate& $\rho v$: Error & Rate& $\rho E$: Error & Rate \\
			\midrule
			16 & 0.0025737 & -- & 0.0045291 & -- & 0.0045291 & -- & 0.040389 & -- \\
			32 & 0.00069247 & 1.89 & 0.0013296 & 1.77 & 0.0013296 & 1.77 & 0.011169 & 1.85 \\
			64 & 0.00019406 & 1.84 & 0.00037977 & 1.81 & 0.00037977 & 1.81 & 0.0029554 & 1.92 \\
			128 & 5.5393e-05 & 1.81 & 0.00011422 & 1.73 & 0.00011422 & 1.73 & 0.00077826 & 1.93 \\
			256 & 1.608e-05 & 1.78 & 3.5496e-05 & 1.69 & 3.5496e-05 & 1.69 & 0.00021095 & 1.88 \\
			512 & 4.7598e-06 & 1.76 & 1.1304e-05 & 1.65 & 1.1304e-05 & 1.65 & 5.8308e-05 & 1.86 \\
			1024 & 1.4667e-06 & 1.70 & 3.709e-06 & 1.61 & 3.709e-06 & 1.61 & 1.6391e-05 & 1.83 \\
			\bottomrule
		\end{tabular}
		\caption{$L^1$ errors of the test in \cref{sec:accuracy} for different resolutions and the corresponding convergence rates. The well-balanced scheme is used with linear reconstruction on the different grids introduced in \cref{sec:slhgrids}.}
		\label{tab:accuracy}
	\end{table}
	
	To test the order of accuracy of the scheme we use a problem from \cite{Xing2013} and \cite{Chandrashekar2017} which involves a known exact solution of the Euler equations with gravity given by
	\begin{align*}
	\rho\args &= 1+\frac{1}{5} \sin(\pi(x+y-t(u_0+v_0))), \qquad u\equiv u_0, \qquad v\equiv v_0,\\
	p\args &= p_0 + t(u_0+v_0)-x-y+\frac{1}{5\pi}\cos(\pi(x+y-t(u_0+v_0))).
	\end{align*}
	The gravitational potential is $\pi(\x)=x+y$, the EoS is the ideal gas EoS \cref{eq:ideal_eos}. In accordance to \cite{Xing2013} and \cite{Chandrashekar2017} we choose $u_0=v_0=1$. $p_0=4.5$. We use our scheme to evolve the initial data with $t=0$ to a final time $t=0.1$ with different grid resolutions on a on all different grid geometries presented in \cref{sec:slhgrids}. $\alpha$ and $\beta$ are set for an isothermal solution (\cref{sec:2d_isothermal}). This is an arbitrary choice, since we want to discuss the accuracy of the scheme if the data are far away from an hydrostatic solution.
	We use MUSCL reconstruction in all tests.
	
	The results are compared to the exact solution at time $t=0.1$. The $L^1$ errors together with the convergence rates are shown in \cref{tab:accuracy}.
	Together with the linear reconstruction our well-balanced scheme is second order accurate on the Cartesian, sinusoidal, and polar grid. On the cubed sphere grid some accuracy is lost due to the grid geometry as the grid is not smooth.

	\section{Conclusions}
	\label{sec:conclusions}
	
	We have proposed a finite volume scheme for the Euler equations with gravitational source term in two spatial dimensions. We have shown that the scheme is well-balanced, and this property holds independently of the EoS, the special hydrostatic solution that has to be maintained, the reconstruction routine, the numerical flux function, and the grid geometry. In numerical tests we could confirm this result by maintaining different hydrostatic solutions with different EoS. We also used different curvilinear meshes and reconstruction routines. Tests over a long simulated time have shown two important features of our scheme: hydrostatic solutions that are stable with respect to convection can be maintained over a long time. Hydrostatic solutions that are unstable with respect to convection develop instabilities, thus reproducing the physical situation. Tests with the cubed sphere grid, which tends to introduce significant discretization errors with non-well-balanced schemes, show the benefit of applying a well-balanced scheme in that case. The well-balanced scheme allows the resolution of small perturbations on a hydrostatic solution. Numerical experiments indicate that our well-balanced scheme is second order accurate if it is combined with a linear reconstruction and a sufficiently accurate time-stepping routine. This holds even on non-Cartesian, non-uniform but smooth curvilinear grids. The generality of the present well-balanced scheme allows us to use the combination of best numerical schemes and grid system that is suitable for a particular problem.
	
	\appendix
	\section{Transformation to curvilinear coordinates}
	\label{appendix:trafo}
	In this section, we will derive the relationship between the two coordinates systems, $(x,y)$ and $(\xi,\eta)$. The differentials are related by
	\[
	\begin{bmatrix}
	\ud x \\
	\ud y \end{bmatrix} = \begin{bmatrix}
	x_\xi & x_\eta \\
	y_\xi & y_\eta \end{bmatrix}
	\begin{bmatrix}
	\ud \xi \\
	\ud \eta \end{bmatrix}, \qquad
	\begin{bmatrix}
	\ud \xi \\
	\ud \eta \end{bmatrix} = \begin{bmatrix}
	\xi_x & \xi_y \\
	\eta_x & \eta_y \end{bmatrix}
	\begin{bmatrix}
	\ud x \\
	\ud y
	\end{bmatrix}
	\]
	and hence
	\[
	\begin{bmatrix}
	\xi_x & \xi_y \\
	\eta_x & \eta_y \end{bmatrix} = \begin{bmatrix}
	x_\xi & x_\eta \\
	y_\xi & y_\eta \end{bmatrix}^{-1} = \frac{1}{J} \begin{bmatrix}
	y_\eta & -x_\eta \\
	-y_\xi & x_\xi \end{bmatrix}, \qquad J = x_\xi y_\eta - x_\eta y_\xi
	\]
	Hence we have proved the following identies
	\begin{equation}
	J \xi_x = y_\eta, \quad J \xi_y = -x_\eta, \quad J \eta_x = -y_\xi, \quad J \eta_y = x_\xi
	\label{eq:trans}
	\end{equation}
	We now transform the conservation law. The flux derivatives transform as
	\[
	\fl_x = \fl_\xi \xi_x + \fl_\eta \eta_x, \qquad \gl_y = \gl_\xi \xi_y + \gl_\eta \eta_y
	\]
	and hence, adding the two equations and making use of~(\cref{eq:trans}), we get
	\begin{eqnarray*}
		\fl_x + \gl_y &=& \frac{1}{J}(y_\eta \fl_\xi - x_\eta \gl_\xi) + \frac{1}{J}(-y_\xi \fl_\eta + x_\xi \gl_\eta) \\
		&=& \frac{1}{J}(y_\eta \fl - x_\eta \gl)_\xi - \frac{1}{J}(y_{\eta\xi}\fl - x_{\eta\xi}\gl) + \frac{1}{J}(-y_\xi \fl + x_\xi \gl)_\eta - \frac{1}{J}(-y_{\xi\eta} \fl + x_{\xi\eta} \gl) \\
		&=& \frac{1}{J}(y_\eta \fl - x_\eta \gl)_\xi + \frac{1}{J}(-y_\xi \fl + x_\xi \gl)_\eta  \qquad \textrm{since} \quad y_{\xi\eta} = y_{\eta\xi}, \quad x_{\xi\eta} = x_{\eta\xi}
	\end{eqnarray*}
	But since $(y_\eta \fl - x_\eta \gl) = (x_\eta^2 + y_\eta^2)^{\tfrac{1}{2}} \hfl$ and $(-y_\xi \fl + x_\xi \gl) = (x_\xi^2 + y_\xi^2)^{\tfrac{1}{2}} \hgl$, and hence we obtain
	\[
	\fl_x + \gl_y = \frac{1}{J} [(x_\eta^2 + y_\eta^2)^{\tfrac{1}{2}} \hfl]_\xi + \frac{1}{J} [ (x_\xi^2 + y_\xi^2)^{\tfrac{1}{2}} \hgl]_\eta
	\]
	This completes the derivation of~(\cref{eq:eul2d}).

	\section*{Acknowledgments}
	This research is supported by the Klaus Tschira Foundation. We would also like to acknowledge the developers of the \emph{Seven-League Hydro Code} (SLH), which has been used for the numerical experiments in this paper with their kind permission. Some of the main developers are Philipp V.~F.~Edelmann from the Heidelberg Institute for Theoretical Studies and Fabian Miczek from the University of W\"urzburg. Praveen Chandrashekar would like to acknowledge the support received from Airbus Foundation Chair on Mathematics of Complex Systems at TIFR-CAM, Bangalore.
	
	\bibliographystyle{abbrv}
	\bibliography{library}
	
\end{document}